\newif\ifpictures
\picturestrue

\newif\ifcomment
\commenttrue

\documentclass[12pt]{amsart}
\usepackage{latex_base}
\usepackage{macros}

\author{Christoph Brauer}

\address{Christoph Brauer, Deutsches Zentrum f\"ur Luft- und Raumfahrt, Institut f\"ur Systemleichtbau, Abteilung Produktionstechnologien, Ottenbecker Damm 12, 21682 Stade, Germany\medskip}

\email{christoph.brauer@dlr.de}

\author{Arne Hindersmann}

\address{Arne Hindersmann, Deutsches Zentrum f\"ur Luft- und Raumfahrt, Institut f\"ur Systemleichtbau, Abteilung Produktionstechnologien, Ottenbecker Damm 12, 21682 Stade, Germany\medskip}

\email{arne.hindersmann@dlr.de}

\author{Timo de Wolff}

\address{Timo de Wolff, Technische Universit\"at Braunschweig, Institut f\"ur Analysis und Algebra, AG Algebra, Universit\"atsplatz 2, 38106 Braunschweig,
 Germany\medskip}

\email{t.de-wolff@tu-braunschweig.de}



\title{Voronoi-Based Vacuum Leakage Detection in Composite Manufacturing}

\begin{document}

\begin{abstract}
In this article, we investigate vacuum leakage detection problems in composite manufacturing.
Our approach uses Voronoi diagrams, a well-known structure in discrete geometry.
The Voronoi diagram of the vacuum connection positions partitions the component surface.
We use this partition to narrow down potential leak locations to a small area, making an efficient manual search feasible.
To further reduce the search area, we propose refined Voronoi diagrams.
We evaluate both variants using a novel dataset consisting of several hundred one- and two-leak positions along with their corresponding flow values.
Our experimental results demonstrate that Voronoi-based predictive models are highly accurate and have the potential to resolve the leakage detection bottleneck in composite manufacturing.
\end{abstract}

\maketitle

\section{Introduction}

Fiber composite materials are widely used in the production of aerospace components and other applications due to their combination of low weight and high stiffness.
Prepreg materials are a specific type of fiber composite that is pre-impregnated with epoxy resin.
They are commonly used to produce large-scale, high-quality structures.
However, the final properties of fiber composite parts depend on more than just the raw materials used.
The production process and its parameters also play a significant role.
To prevent voids and ensure uniform consolidation, the fiber preform is sealed in an airtight cavity that is evacuated during the curing process.
This cavity usually consists of a rigid mold and a flexible bag, which is sealed to the mold edges with tacky tape.
A variety of auxiliary materials in the cavity, such as a breather layer, ensure even evacuation by the vacuum pump.
For large components, several vacuum connections are typically required on the mold side to maintain uniform pressure distribution across the entire surface.   
\medskip

In this study, we examine the prepreg manufacturing process as an example.
However, we see no obstruction to apply the proposed leakage detection methodology also to the vacuum infusion manufacturing process.
In the prepreg process, the component is cured in an autoclave at elevated temperature and pressure after being sealed.
If the vacuum bag, tacky tape, mold, or vacuum connection leaks during this process, voids are likely to form.
Such voids can permanently impair the mechanical performance of the composite structure and are therefore highly undesirable.
To prevent this, the cavity is tested for leaks before curing.
In principle, either pressure or volume flow measurements can be used for this purpose.
In the \struc{pressure measurement} method, the connection between the vacuum pump and the cavity is closed, and the subsequent pressure increase is evaluated.
If the pressure rise over time is small, the cavity can be considered airtight. \cite{haschenburger2019a, dlr190099}
For \struc{volume flow measurement}, a flow sensor is integrated into the vacuum line to determine the volume flow rate.
If the measured flow is very low or zero, the cavity can be considered airtight.\cite{dlr190099, haschenburger2021a}
\medskip

While both pressure and volume flow measurements provide an indication of the cavity’s tightness, neither can identify the precise location of a leak.
Therefore, additional methods must be used to detect leaks. The following approaches are well established in the literature and are also applied in industrial environments.
The \struc{helium leak test} uses a chamber that is placed on the vacuum bag and filled with helium.
A sensor measuring the helium concentration is installed upstream of the vacuum pump.
If a high helium concentration is detected, a leak is likely present in the area covered by the chamber.
The chamber is then moved to the next section until the entire component has been inspected. \cite{haschenburger2019a}
The \struc{ultrasonic leak test} exploits the fact that small leaks often generate high-frequency noise due to incoming air.
A directional microphone is used to convert these inaudible frequencies into lower, audible ones.
Since the microphone is portable, even large components can be examined.
The \struc{thermographic leak test} is based on the Joule-Thomson effect: leaks cause localized cooling of the surrounding air.
An infrared camera visualizes these cooled areas and thus identifies the leak locations.\cite{DLRK2011}
\medskip

All of these methods allow for efficient localization of leaks in relatively small or geometrically simple components.
However, larger components, or those with complex geometrical features, typically require significantly more time and effort to inspect.
As a result, leakage detection can become a bottleneck in the production of large structures such as wing covers.
One strategy to address this challenge is to employ an additional technique that narrows the search area on large components to a size that can be inspected more quickly.
In the past, separate flow measurements at different vacuum connections were used for this purpose.
The underlying intuition is that air follows the \struc{path of least resistance} from a leak location to the nearest vacuum connection.
Consequently, most of the airflow will take this path, and the flow meter at the corresponding connection will indicate the highest value.
In other words, if leaks are present, at least one of them is likely to be located closest to the vacuum connection with the largest measured flow.
This concept can be formalized using a well-established, classical object from discrete geometry called \struc{Voronoi diagram}.
\medskip

Previous literature on methodologies for transforming volumetric flow measurements into leak localization predictions includes knowledge-based approaches as well as machine learning models.
In \cite{haschenburger2021a}, the authors investigate a hand-crafted mapping between order relation patterns in three flow measurements and specific subregions of a triangular surface.
The studies \cite{dlr190099, brauer2022a, haschenburger2022computational} rely on the same dataset captured on a 1.5-square-meter, quadratic experimental part with four vacuum connections.
In \cite{dlr190099}, the authors subdivide the quadratic surface into regular grids of varying granularity and use machine learning-based \struc{classification} models to predict whether each grid cell contains a leak.
This principally enables the localization of multiple leaks simultaneously.
The authors of \cite{brauer2022a} integrate symmetry information into \struc{regression} neural networks to predict single-leak coordinates.
The comparative study \cite{haschenburger2022computational} considers potential flow theory, different regression models and a volumetric flow matching method, all for the task of single-leak detection.
The master's thesis \cite{naveenachandran2023}, which was supervised by the authors of this work, investigates regression neural networks for multi-leak detection using the same dataset as in this article.
\medskip

In this paper, we propose leakage detection models based on Voronoi diagrams.
Our work generalizes the partitioning approach from \cite{haschenburger2021a}, which uses a Voronoi diagram on a triangular surface without explicitly establishing this link, to arbitrary two-dimensional surfaces.
Furthermore, we introduce \struc{refined Voronoi diagrams}, a new generalization of Voronoi diagrams that allows for a finer surface partition.
Similar to \cite{brauer2022a} our approach incorporates prior information about the part geometry without being limited to a specific geometry, such as a square.
Technically, our predictive models are classification models.
However, unlike \cite{dlr190099}, they are purely knowledge- and geometry-based, and, instead of using a regular grid, we use a Voronoi-based subdivision that exploits the geometric features of the part and the vacuum setup.
Our numerical experiments, conducted on an industrial-scale wing dataset collected at the DLR Institute of Lightweight Systems in Stade, Germany, demonstrate that Voronoi diagrams are well-suited for leakage localization and that our methods predict Voronoi cells containing leaks with high confidence.
We provide our wing dataset and the code necessary to reproduce the subsequently reported results at 
\begin{quote}
	\url{https://github.com/chrbraue/leakage\_detection}.
\end{quote}
\medskip

The remainder of this paper is structured as follows.
In Section~\ref{sec:methodology}, we present our mathematical model of the vacuum setup and introduce classical and refined Voronoi diagrams, as well as respective predictive models for leak localization.
Section~\ref{sec:experiments} includes the experimental evaluation of the proposed methods.
In Section~\ref{sec:discussion}, we discuss and interpret the obtained results with a focus on practical implications.
Section~\ref{sec:conclusion} concludes the main part of the paper.
Subsequently, Appendix~\ref{appendix:discrete_geometry} introduces basic concepts from discrete geometry, while Appendix~\ref{appendix:voronoi_diagrams} provides additional information on classical Voronoi diagrams.
Appendix~\ref{appendix:refined_voronoi_diagrams} offers a broader definition of refined Voronoi diagrams and an algorithm for constructing them since, to the best of our knowledge, they are a novel structure in both leak detection and discrete geometry (even though it is closely related to higher-order Voronoi diagrams).
\medskip

\subsection*{Notation}

The following notation is used throughout the paper.
The surface of the vacuumized part is represented by a connected two-dimensional point set $\struc{\surface}$, and $\struc{\pumpset}$ denotes a discrete two-dimensional point set representing the positions of vacuum connections.
Vectors (namely $\pb_i \in \pumpset$, flow measurements $\xb \in \R^k$, and leak positions $\yb\in\R^2$) are printed in bold, and scalars are in normal font.
We denote covers of $\surface$ by $\struc{\voronoidiagram} := \{\voronoicell_1,\dots, \voronoicell_c\}$ (or similar), where the elements $\voronoicell_1,\dots, \voronoicell_c \subset \surface$ are disjoint subsets that together exactly cover $\surface$.
To denote functions that map a vector of flow values to an element of a cover, we use $f: \R^k \rightarrow \voronoidiagram$.
Furthermore, we abbreviate $\struc{[k]} := \{1,\ldots,k\}$ and use $[k] \times [k]$ to denote the Cartesian product of the same set with itself. Finally, $\struc{\Vert\cdot\Vert}$ refers to the Euclidean norm.
Terms are introduced in great detail below as needed.

\section{Methodology}
\label{sec:methodology}

In this section, we present our mathematical model of the vacuum setup.
The model incorporates the surface of the evacuated component, the leak location, the locations of the vacuum connections, and the sensor data collected at these connections.
We further introduce predictive models that reduce the search area for leaks from the entire surface to smaller subregions.
We provide only the mathematical details necessary to understand the proposed methods.
The appendix contains rigorous definitions and derivations of the underlying mathematical structures.

\subsection{Mathematical model of the vacuum setup}

The methods presented in this work focus on flat components or components with negligible curvature.
Accordingly, the surface of the evacuated component is modeled as a two-dimensional set $\struc{\surface} \subset \R^2$.
If the actual surface is not perfectly flat, $\surface$ denotes its two-dimensional projection.
The positions of the vacuum connections are represented as a finite set of points
\begin{equation*}
    \struc{\pumpset := \left\{\pb_1,\dots,\pb_k \right\}} \subset \surface
\end{equation*}
and a leak in the vacuum bag is represented by a two-dimensional coordinate vector $\struc{\yb = (y_1, y_2)}\in\surface$.
Since leaks are typically very small, their area or diameter is neglected in this model.
Alternatively, $\yb$ may be interpreted as the geometric center of a leak.
\medskip

To describe the relationship between a leak and the corresponding flow measurements, we consider a leak located at position $\yb \in \surface$.
This leak causes characteristic flow values $\struc{\xb = (x_1,\dots,x_k)}$ that are measured at the vacuum connections located at $\pb_1,\dots,\pb_k$.
Although measurement noise and minor variations between different realizations of the same vacuum setup may cause small deviations in the measured values, we do not explicitly take these into account.

\subsection{Predictive models for leak localization}
\label{sec:predictive_models_for_leak_localization}

We present two approaches for narrowing the search area for potential leaks from the entire component surface to smaller, more manageable regions.
This idea is motivated by both practical experience and mathematical considerations.
As discussed above, methods such as helium, ultrasonic, and thermographic leak testing can efficiently locate leaks in small areas, but are less effective for large or complex components such as wing covers.
From a mathematical perspective, developing a model that predicts exact leak coordinates is inherently difficult, even in the case of a single leak.
Instead, it is often more meaningful to determine a confidence region, which effectively narrows down the search area.
When multiple leaks occur simultaneously, as is often the case in practice, the problem is even more complex.
The same set of flow measurements can then potentially be caused by different numbers or distributions of leaks.
Under appropriate assumptions, this ambiguity can be proven rigorously, but a detailed analysis is beyond the scope of this paper.
\medskip

Let $\struc{\voronoidiagram := \{\voronoicell_1,\dots,\voronoicell_c\}}$ denote a cover of the surface $\surface$, that is,
\begin{equation*}
    \bigcup_{i=1,\dots,c}\voronoicell_i = F \quad \text{and} \quad \voronoicell_i \cap \voronoicell_j = \emptyset \ \text{for} \  i\neq j\, .
\end{equation*}
In Sections~\ref{sec:classic_voronoi_diagrams} and \ref{sec:refined_voronoi_diagrams}, we introduce two concepts for constructing such a cover based on Voronoi diagrams, together with functions $\struc{\backwardvoronoi}: \R^k \to \voronoidiagram$ that map a measurement vector $\xb$ to a subset $\backwardvoronoi(\xb) = \voronoicell_i$ of the surface that is predicted to contain a leak.
The choice of the cover $\voronoidiagram$ is purely geometric and knowledge driven; it does not require any observed data.
In Section~\ref{sec:strategies_for_multi_leak_localization}, we present strategies that generalize the two concepts to the multi-leak case.

\subsubsection{Classic Voronoi diagrams}
\label{sec:classic_voronoi_diagrams}

\begin{figure}[t]
  \centering
  \includegraphics[width=.45\textwidth, trim=12 12 12 12, clip]{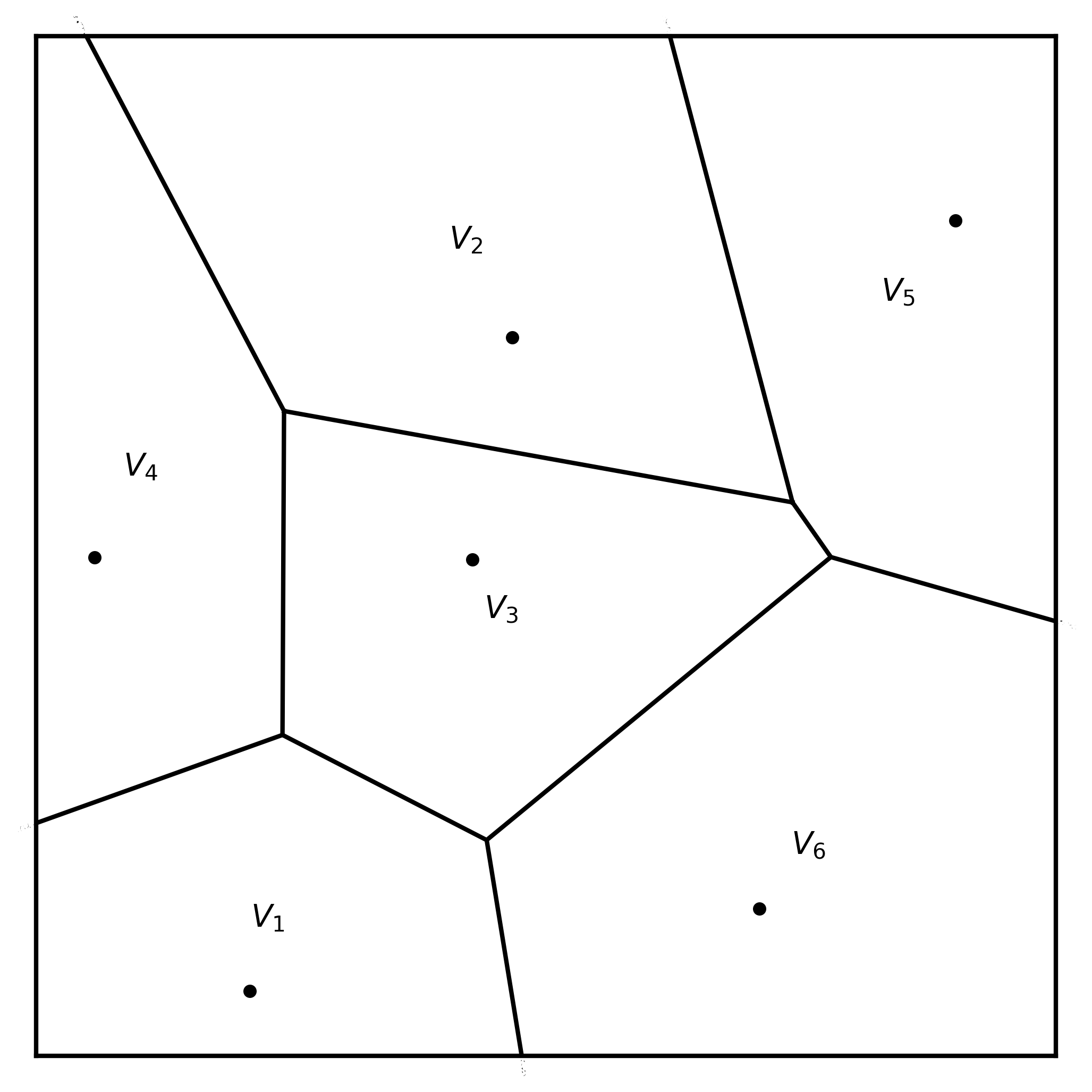}
  \caption{Classic Voronoi diagram of a random set of points}
  \label{fig:classic_voronoi_diagram_random_example_main_text}
\end{figure}

Voronoi diagrams, introduced by Georgy Voronoi in 1908 (and, in the special cases of two and three dimensions, already by Dirichlet in 1850), are a well-known construct in discrete geometry; see e.g. \cite{aurenhammer2013voronoi}.
Voronoi diagrams have numerous applications in fields such as meteorology (e.g., famously applied by Alfred H. Thiessen for the prediction of rainfall already in 1911; also referred to as \struc{Thiessen polygons}; see, e.g., \cite{Brassel:ThiessenPolygons}) biology (e.g., \cite{Indermitte:VoronoiBio}), and crystallography  (e.g., \cite{Blatov:VoronoiCrystal}), to name a few.
Below, we present the standard definition of a two-dimensional Voronoi diagram, as it is commonly found in the literature.
For a more in-depth introduction to discrete geometry and Voronoi diagrams, refer to Appendices \ref{appendix:discrete_geometry} and \ref{appendix:voronoi_diagrams}, respectively.
Since the following section introduces a refinement of Voronoi diagrams, we will sometimes refer to the version defined here as the \struc{classic Voronoi diagram}.
\medskip

\begin{definition}
	Let $\pumpset := \{\pb_{1},\dots,\pb_{k}\} \subset \R^2$ be a finite set of points.
    For each $i \in [k]$, the corresponding \struc{Voronoi cell} is defined as
	\begin{equation*}
		\struc{\voronoicell_{i}} \coloneqq \left\{\yb\in\R^{2} \;\middle|\; ||\yb - \pb_{i}|| \leq ||\yb - \pb_{j}|| \text{ for all } j\neq i\right\} \, .
	\end{equation*}
	The collection
    \begin{equation*}
        \struc{\voronoidiagram(\pumpset)} \coloneqq \{\voronoicell_1,\dots,\voronoicell_k\}
    \end{equation*}
    is called the \struc{Voronoi diagram} of $\pumpset$.
    \label{def:voronoi_diagram_2d}
\end{definition}
\medskip

\begin{example}
    Figure~\ref{fig:classic_voronoi_diagram_random_example_main_text} illustrates the classic Voronoi diagram of a random set of points $\pumpset = \{\pb_{1},\dots,\pb_{6}\}$ sampled from a uniform distribution on $\surface = [-1, 1]^2$.
    As can be seen, the Voronoi diagram $\voronoidiagram(\pumpset) = \{\voronoicell_1,\dots,\voronoicell_6\}$ partitions the surface into six subsets.
    The line segment between two Voronoi cells is the set of points equidistant from the two respective points from $\pumpset$.
    Accordingly, the nodes at which three line segments meet are equidistant from three points.
    \label{example:voronoi_diagram}
\end{example}
\medskip

We now introduce our first predictive model based on classic Voronoi diagrams.
This model relies on the assumption that, if a leak is present, the largest flow will be measured at the vacuum connection closest to the leak.
Conversely, when flow measurements are available and the leak location is to be inferred, the leak is expected to lie within the region of points that are closest to the vacuum connection exhibiting the largest flow value.
Hence, the Voronoi diagram constructed from the vacuum connection locations provides a natural geometric framework to formalize and implement this search rule.
If $x_i$ is the largest entry in the flow vector $\xb$, then the corresponding Voronoi cell $\voronoicell_i$ is predicted to contain the leak.
\medskip

\begin{definition}
Let $\voronoidiagram(\pumpset) = \{\voronoicell_1, \dots, \voronoicell_k\}$ denote the classic Voronoi diagram of the vacuum connection locations $\pb_1, \dots, \pb_k$.
Define
\begin{equation*}
    \struc{i^*(\xb)}
\end{equation*}
as the index corresponding to the largest entry of the measurement vector $\xb \in \R^k$:
\begin{equation*}
    x_{i^*(\xb)} \geq x_j \quad \text{for all}\quad j \in [k] \setminus \{i^*(\xb)\}.
\end{equation*}
Then,
\begin{equation*}
    \backwardvoronoi_1(\xb) \coloneqq \voronoicell_{\,i^*(\xb)}
\end{equation*}
defines a function $\backwardvoronoi_1 : \R^k \to \voronoidiagram(\pumpset)$, referred to as the \struc{classic Voronoi predictor}, 
which maps a measurement vector $\xb \in \R^k$ to the Voronoi cell associated with the vacuum connection exhibiting the largest flow measurement.
\end{definition}
\medskip

\subsubsection{Refined Voronoi diagrams}
\label{sec:refined_voronoi_diagrams}

Various generalizations of the classic Voronoi diagram have been proposed and applied in the literature.
These include generalizations with respect to the underlying metric (other than the Euclidean $\ell_2$-norm), cells defined relative to objects other than single points (e.g., line segments), and extensions to spatial dimensions beyond two; see, for example, \cite{aurenhammer2013voronoi} for an overview.
Among these, the most extensively studied variant is arguably the \struc{higher-order Voronoi diagram}, also referred to as the \struc{$k$-nearest neighbor Voronoi diagram}.
These objects were first introduced in \cite{shamos1975closest}, and remain an active field of research; see also, e.g., \cite{ChazelleEdelsbrunner,Chan2024OptimalHOVD}.
Order-$d$ Voronoi cells are defined such that each point within a given cell shares the same $d$ nearest sites, irrespective of their order.
\medskip

In this work, we introduce \struc{refined Voronoi diagrams}, which, to the best of our knowledge, constitute a novel generalization and a refinement of both classic and higher-order Voronoi diagrams.
The key distinction is that, in refined Voronoi diagrams, the ordering of the nearest sites is explicitly taken into account -- each cell corresponds to a particular permutation of the $d$ nearest sites.
In the following, we present the definition of order-two refined Voronoi diagrams, which are also employed in our empirical experiments.
A more general definition for arbitrary order-$d$ refined Voronoi diagrams, along with details on their construction, is provided in Appendix~\ref{appendix:refined_voronoi_diagrams}.
\medskip

\begin{definition}
    Let $\pumpset = \{\pb_{1},\dots,\pb_{k}\} \subset \R^2$ be a finite set of points and
    \begin{equation*}
        \struc{\sequence_2}\coloneqq \left\{(i_1, i_2)\in [k]\times [k] \;\middle|\; i_1 \neq i_2\right\}
    \end{equation*}
    the set of ordered index pairs without repetition.
    For each $\struc{\tb} = (i_1, i_2)\in \sequence_2$, the corresponding \struc{order-two refined Voronoi cell} is defined as
    \begin{equation*}
        \struc{\voronoicell_{\tb}} \coloneqq \left\{\yb \in \R^2 \;\middle|\; \Vert\yb - \pb_{i_{1}}\Vert\leq \Vert\yb - \pb_{i_{2}}\Vert\leq \Vert\yb - \pb_{j}\Vert \text{ for all } j \in [k] \setminus \{i_1, i_2\}\right\} .
    \end{equation*}
    The collection
    \begin{equation*}
        \struc{\generalizedvoronoidiagram_2(\pumpset)} \coloneqq \{\voronoicell_{\tb}\}_{\tb \in \sequence_{2}}
    \end{equation*}
    is called the \struc{order-two refined Voronoi diagram} of $\pumpset$.
\end{definition}
\medskip

\begin{figure}[t]
  \centering
  \includegraphics[width=.45\textwidth, trim=12 12 12 12, clip]{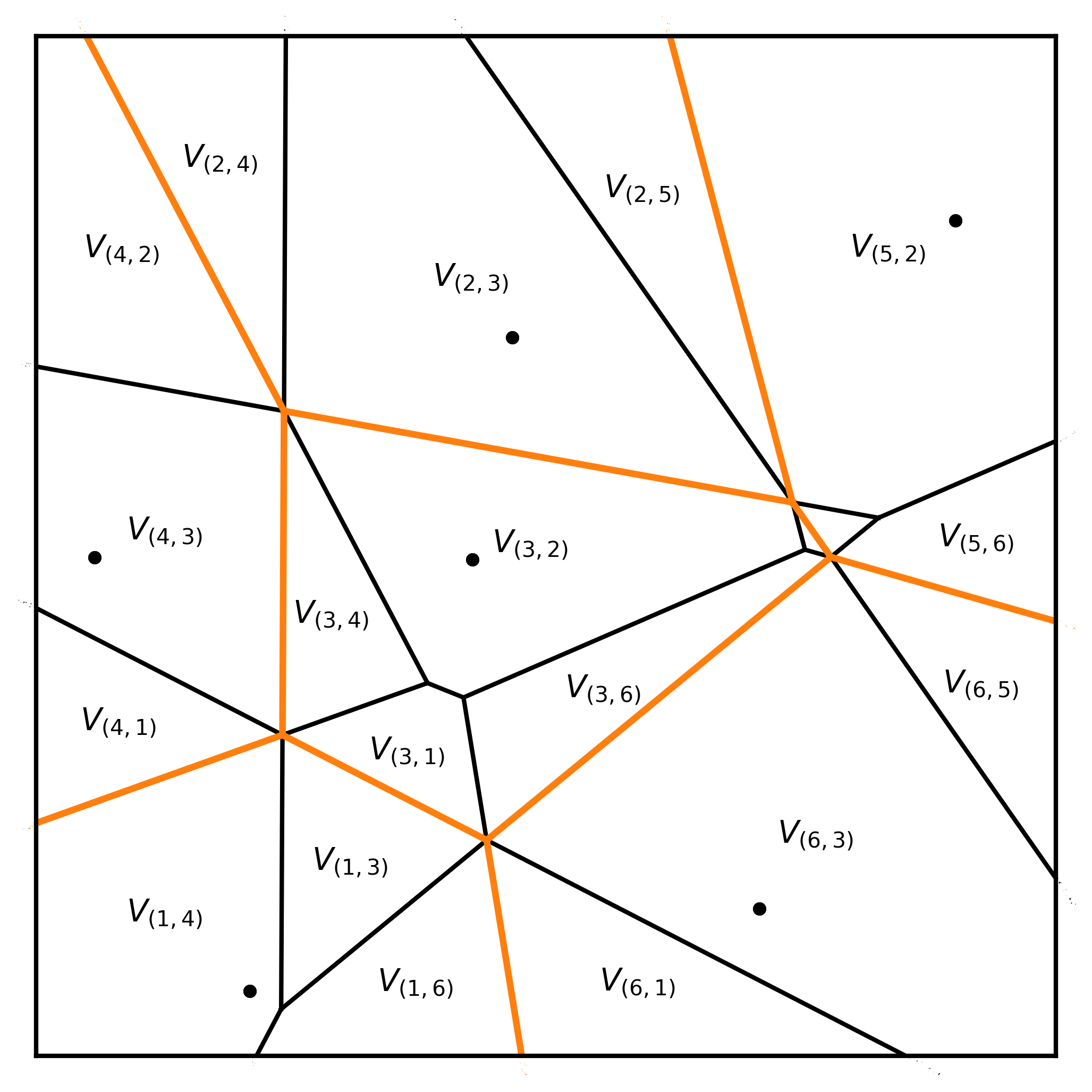}
  \caption{Order-two refined Voronoi diagram of a random set of points}
  \label{fig:refined_voronoi_diagram_random_example_main_text}
\end{figure}

\begin{example}
Figure~\ref{fig:refined_voronoi_diagram_random_example_main_text} revisits the example introduced in Figure~\ref{fig:classic_voronoi_diagram_random_example_main_text}.
The classic Voronoi diagram is shown in orange.
Consider the central cell $\voronoicell_3$, which has five edges.
In the refined Voronoi diagram (black), $\voronoicell_3$ is subdivided into five refined Voronoi cells.
Each of these subcells, denoted $\voronoicell_{(3, j)}$, corresponds to one of the five neighboring cells and contains the points that are closest to $\pb_3$ and second-closest to $\pb_j$.
While $V_3$ is adjacent to each other classic Voronoi cell, this is not true for the remaining five cells.
Each of them has only three adjacent cells and, accordingly, three edges; thus, each is subdivided into three refined Voronoi cells.
For instance, the cell $\voronoicell_1$ is subdivided into $\voronoicell_{(1, 3)}$, $\voronoicell_{(1, 4)}$ and $\voronoicell_{(1, 6)}$, whereas $\voronoicell_{(1, 2)}$ and $\voronoicell_{(1, 5)}$ are empty.
The connection to the corresponding higher-order (2-nearest neighbor) Voronoi diagram is as follows: each pair of refined Voronoi cells $\voronoicell_{(i, j)}$ and $\voronoicell_{(j, i)}$ is separated by an edge that also belongs to the classic Voronoi diagram.
The corresponding higher-order Voronoi diagram can then be obtained by merging each such pair of cells.
Informally, this means that removing all orange edges in Figure~\ref{fig:refined_voronoi_diagram_random_example_main_text} yields the corresponding higher-order (2-nearest neighbor) Voronoi diagram.
\end{example}
\medskip

With the definition of a refined Voronoi diagram in hand, introducing our second predictive model is straightforward.
Building on our earlier assumption that a leak is closest to the vacuum connection with the largest flow, we extend this idea as follows:
We expect the leak to be closest to the vacuum connection with the largest flow and second-closest to the connection with the second-largest flow.
This notion of first- and second-closest proximity is precisely captured by order-two refined Voronoi diagrams.
The following definition adapts the concept of the Voronoi predictor from classic Voronoi diagrams to order-two refined Voronoi diagrams.
\medskip

\begin{definition}
\label{def:refined_voronoi_predictor}
Let $\generalizedvoronoidiagram_2(\pumpset) \coloneqq \{\voronoicell_{\tb}\}_{\tb \in \sequence_{2}}$ denote the order-two refined Voronoi diagram of the vacuum connection locations $\pb_1, \dots, \pb_k$.
Define
\begin{equation*}
    \struc{\tb^*(\xb)} \coloneqq (i^*_1(\xb), i^*_2(\xb))
\end{equation*}
as the indices corresponding to the two largest entries of $\xb$, sorted in descending order:
\begin{equation*}
x_{i^*_1(\xb)} \ge x_{i^*_2(\xb)} \ge x_j \quad \text{for all } j \in [k] \setminus \{i^*_1(\xb), i^*_2(\xb)\}.
\end{equation*}
Then,
\begin{equation*}
    \backwardvoronoi_2(\xb) \coloneqq \voronoicell_{\,\tb^*(\xb)}
\end{equation*}
defines a function $\backwardvoronoi_2 : \R^k \to \generalizedvoronoidiagram_2(\pumpset)$, referred to as the \struc{order-two refined Voronoi predictor}, 
which maps a measurement vector $\xb \in \R^k$ to the order-two refined Voronoi cell associated with the vacuum connections exhibiting the two largest flow measurements.
\end{definition}
\medskip

\subsubsection{Strategies for multi-leak localization}
\label{sec:strategies_for_multi_leak_localization}

Both the Voronoi and refined Voronoi predictors map a set of flow values to a single region predicted to contain a leak.
In practice, however, the vacuum bag often contains more than one leak, and the total number of leaks is unknown.
To apply our predictive models to such situations, we introduce the strategy of \struc{repeatedly applying a Voronoi predictor}.
If more than one leak is possible, this strategy assumes that the region predicted by the Voronoi predictor contains at least one of them.
Once the first leak is located and repaired, flow measurements change until a new flow equilibrium is established.
Then, we use another application of the Voronoi predictor to locate the next leak.
This procedure is repeated until the total flow drops below a certain threshold, which indicates that the setup is sufficiently tight.
\medskip

While the above-mentioned strategy involves fixing a leak before predicting the next one based on another measurement, it is, of course, desirable to use a single initial measurement to \struc{detect multiple leaks simultaneously}.
We can, for example, predict the classic Voronoi cells corresponding to the two connections with the largest flows to detect two leaks at once with a single measurement.
This principle is similar to that of the refined Voronoi predictor, except the prediction consists of two classic Voronoi cells instead of one refined Voronoi cell.
However, this strategy's capability is limited if the number of leaks is unknown.
Furthermore, it is not equally applicable to the classic and refined Voronoi predictors.
Therefore, we do not consider it in our experimental evaluation.
Another way to predict multiple leaks based on a single measurement is training a deep learning \cite{goodfellow2016} model.
Multi-label classification models, as used in combination with regular grid cells in \cite{dlr190099}, are principally suitable for this purpose.
Each grid cell -- or, in our case, each classic or refined Voronoi cell -- then corresponds to one class.
We will discuss the applicability of these models further below.
However, implementing machine learning models for our current use case is beyond the scope of this paper and will be the subject of future work.

\section{Experiments}
\label{sec:experiments}

\begin{figure}[t]
\centering

\begin{subfigure}{\linewidth}
	\centering
	\includegraphics[width=.9375\linewidth]{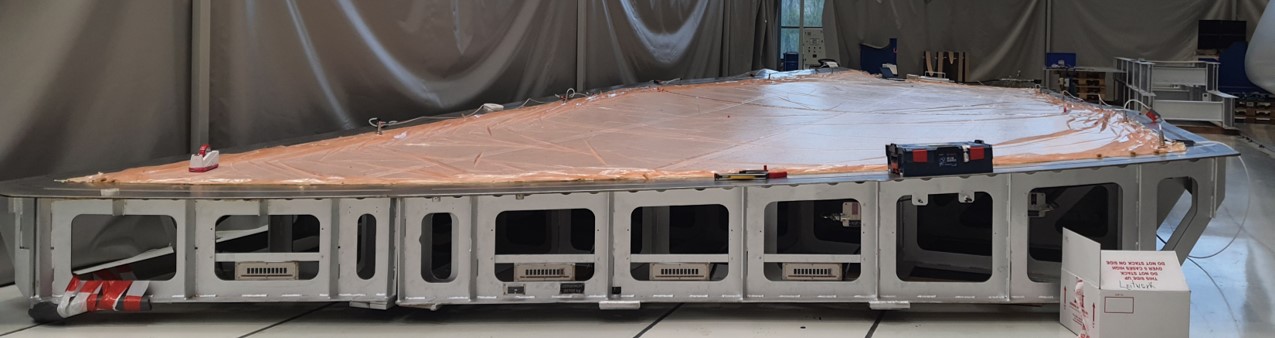}\\
	\caption{Upper wing cover mold}
	\label{fig:layup_tool}
\end{subfigure}

\vspace{1em}

\begin{subfigure}{\linewidth}
	\centering
	\includegraphics[width=.9375\linewidth]{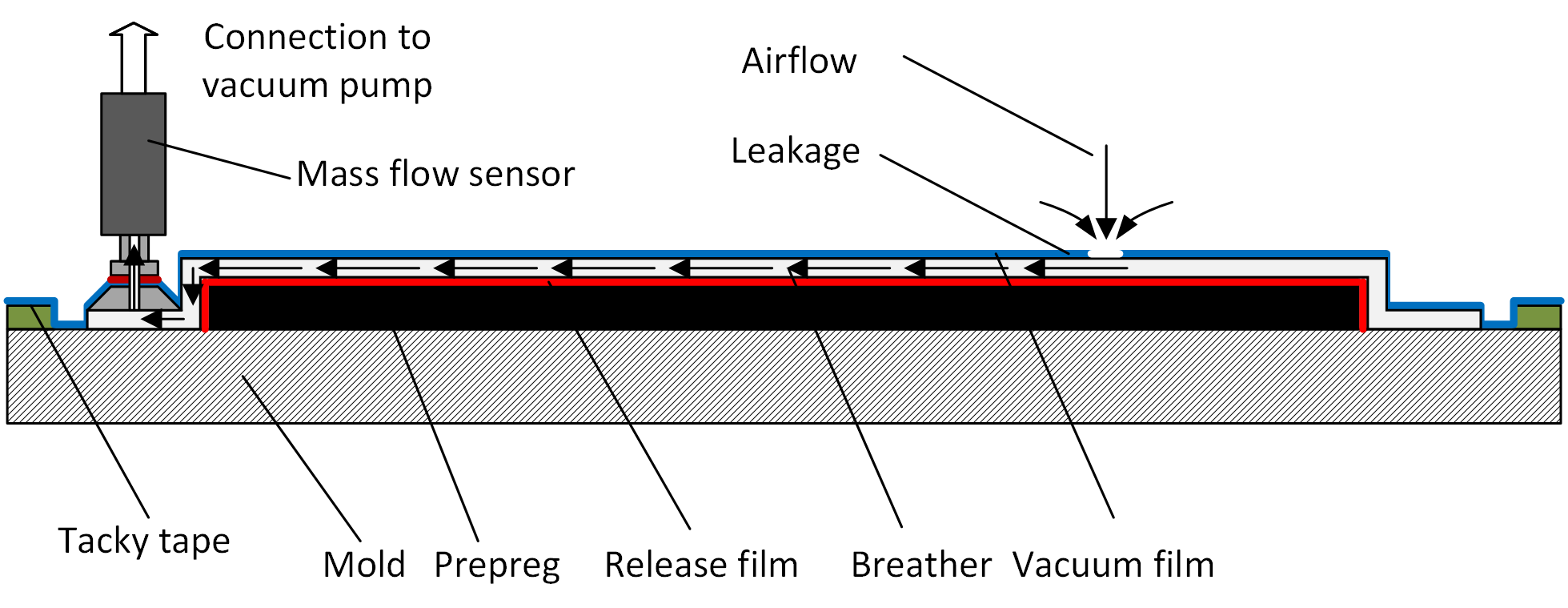}\\
	\caption{Schematic cross-sectional view}
	\label{fig:ExpSetUp}
\end{subfigure}

\caption{Experimental setup}
\label{fig:experimental_setup}

\end{figure}

We evaluate the predictive models introduced in the previous section using the upper wing cover mold shown in Figure~\ref{fig:layup_tool}.
This mold, manufactured for research purposes, measures approximately $5.2\times 16$ meters.
It is industrial in scale, except for the missing wingtip at the narrower end. 
Before data acquisition, we first cover the tool surface with breather material (Airtech Ultraweave 1332), followed by vacuum film (Airtech WL 600V).
Then, we seal the setup airtight along the boundary with tacky tape (Airtech GS 213-3).
Finally, we install ten vacuum connections (Airtech Multi-Valve-407) equipped with flow sensors (Vögtlin Instruments GmbH red-y smart GSM) at specified positions around the perimeter.
Figure~\ref{fig:ExpSetUp} shows a schematic cross-sectional view of the experimental setup.
We assume that the air flow caused by a leak occurs only within the breather.
This assumption is based on the fact that the release film (red) between the breather and the prepreg (black) acts as a barrier.
For this reason, we omit the prepreg material from the actual experimental setup.
\medskip

\subsection{Data acquisition}

\begin{figure}[t]
  \centering
  \includegraphics[height=.2875\linewidth]{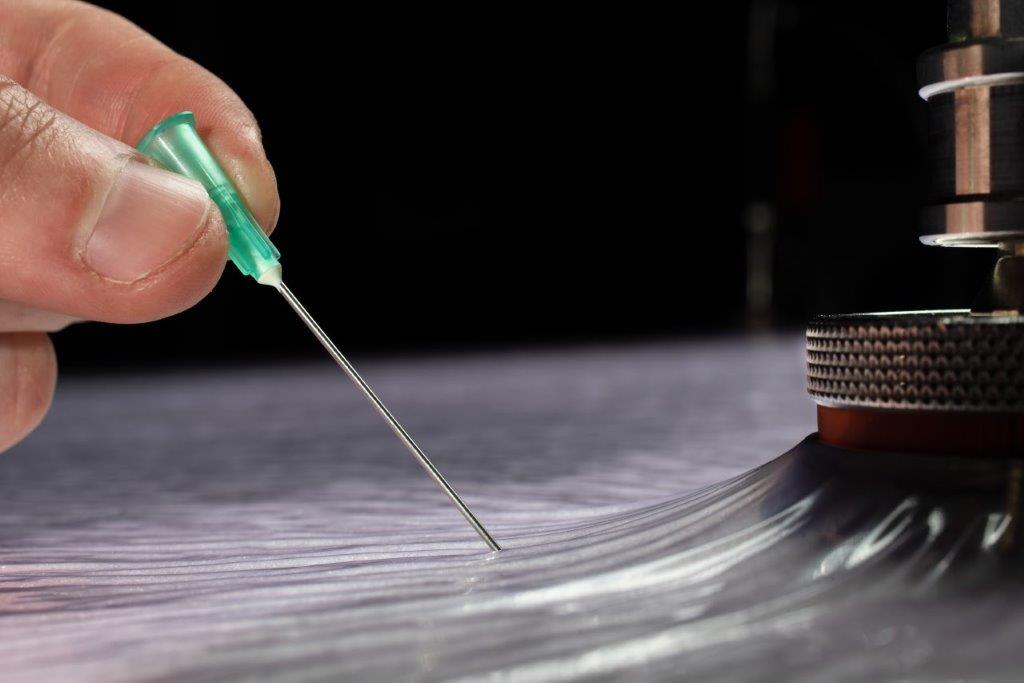}\quad
  \includegraphics[height=.2875\linewidth]{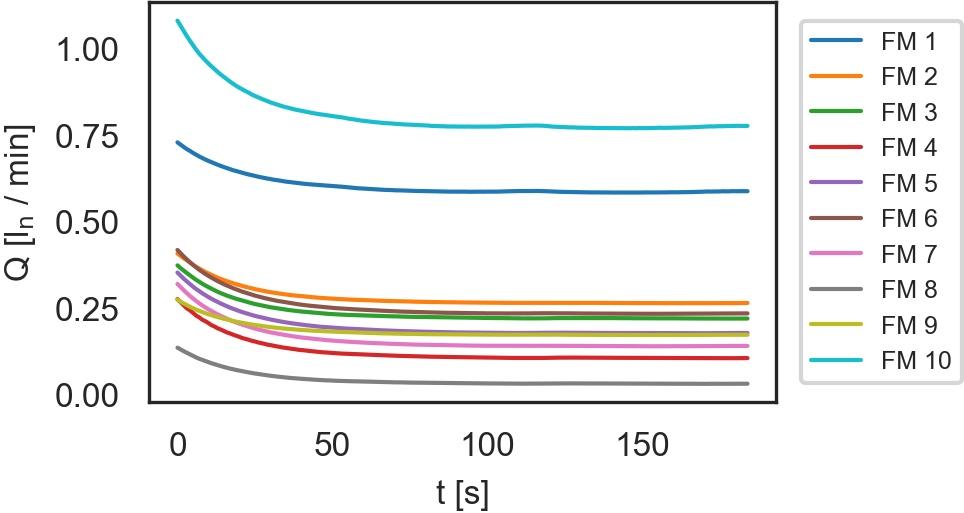}
  \caption{Leak introduction and measurements (FM = Flow Meter)}
  \label{fig:leak_introduction}
\end{figure}

Once the setup is airtight and evacuated, the main process to acquire a measurement is as follows:
While the vacuum pump is running, we use a hypodermic needle to pierce the vacuum bag at a specific position, as illustrated in the left part of Figure~\ref{fig:leak_introduction}.
Air begins to flow through the vacuum connections due to the induced leak.
After a short period, the airflow at all connections stabilizes.
We record the equilibrium values (see the right part of Figure~\ref{fig:leak_introduction}) along with the leak position.
Finally, we seal the leak with adhesive tape and repeat the process at the next location.
\medskip

\begin{figure}[t]
  \centering
  \includegraphics[trim=5 5 5 5, clip]{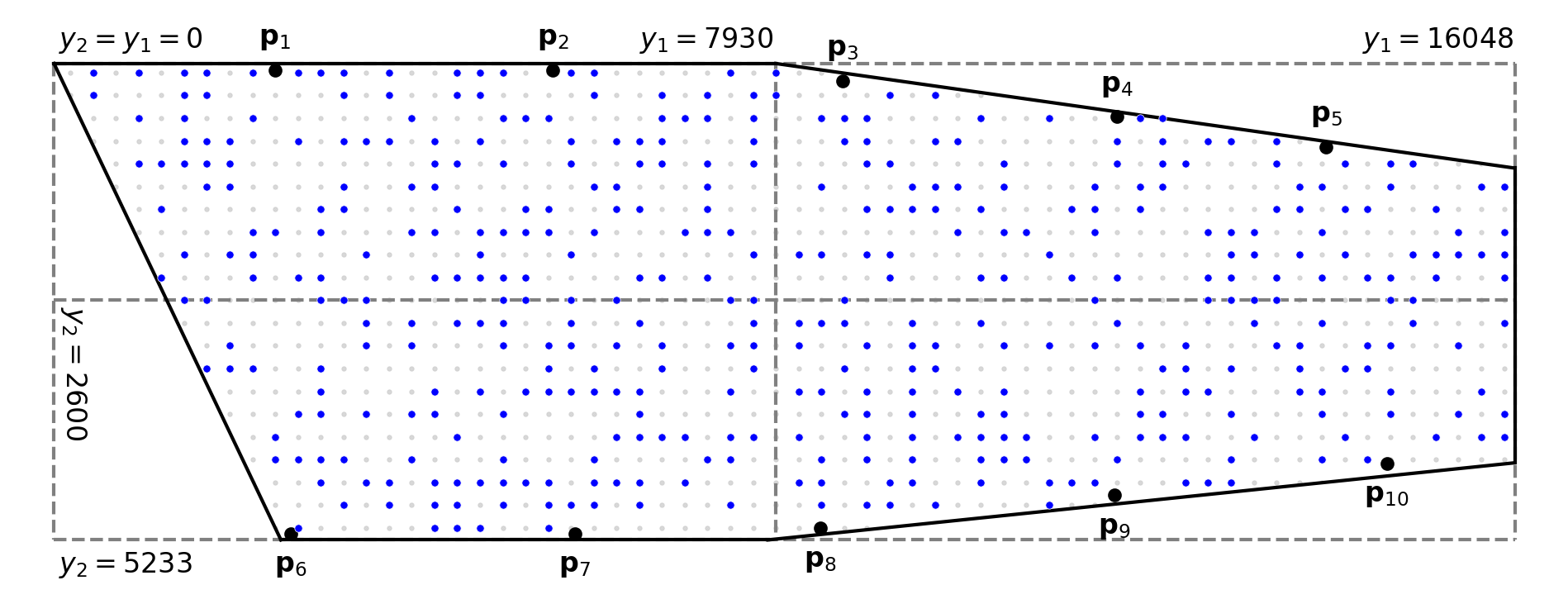}
  \caption{Experimental data}
  \label{fig:experimental_data}
\end{figure}

Figure~\ref{fig:experimental_data} shows the detailed positions of the measurement points.
To generate leak locations, we use a regular grid with a step size of $250$ mm, and select a random grid node for each trial.
To efficiently produce two-leak samples, we take an intermediate step between two single-leak trials.
Instead of immediately fixing a previously recorded single leak, we first introduce a second leak at the next selected grid node.
Then, we record the equilibrium airflow along with the positions of both open leaks.
After this step, we fix the previously recorded single leak, and measure the equilibrium airflow corresponding to the second single leak.
We repeat this procedure as needed.
Using this approach, we generated a total of $413$ single-leak examples and $440$ two-leak examples (some positions appear in more than one two-leak pair).

\subsection{Metrics}

We evaluate our predictive models in terms of two metrics: \struc{accuracy} and \struc{mean Euclidean distance}.
Accuracy is defined as the proportion of correct predictions among all evaluated samples. Formally, this is expressed as follows:

\begin{equation*}
    \mathrm{Accuracy} = \frac{\text{Number of Correct Predictions}}{\text{Number of Samples}}.
\end{equation*}
\smallskip

The Euclidean distance is commonly used to measure the discrepancy between two vectors, e.g., between true and predicted coordinates.
However, since we are considering categorical predictions, in which each Voronoi cell corresponds to one category, the use of the Euclidean distance requires further explanation.
Ultimately, we evaluate Euclidean distances between true leak coordinates and their projections onto predicted Voronoi cells.
The \struc{Euclidean projection} of a point $\yb$ onto a Voronoi cell $\voronoicell_i$ is defined as

\begin{equation*}
    \mathrm{proj}_{\voronoicell_i}(\yb) = \arg\min_{\zb\in \voronoicell_i} \Vert \zb - \yb \Vert \, .
\end{equation*}
\smallskip

\noindent In other words, the Euclidean projection of $\yb$ onto $\voronoicell_i$ is the point in $\voronoicell_i$ that is closest to $\yb$ in terms of Euclidean distance.
Based on this, we use the following definition in our evaluation:

\begin{equation*}
    \text{Mean Euclidean Distance} = \frac{1}{\text{Number of Samples}} \sum_{(\xb, \yb)\in \text{Samples}} \Vert \yb - \mathrm{proj}_{\backwardvoronoi(\xb)}(\yb) \Vert
\end{equation*}
\smallskip

\noindent This is precisely the mean distance between the actual leak coordinates and their projections onto the predicted Voronoi cells, where either $\backwardvoronoi=\backwardvoronoi_1$ (classic Voronoi predictor) or $\backwardvoronoi=\backwardvoronoi_2$ (refined Voronoi predictor).
For each sample $(\xb, \yb)$, there are two possible cases.
First, if the Voronoi predictor makes a correct prediction, then the actual leak position is inside the predicted Voronoi cell.
In this case, $\mathrm{proj}_{\backwardvoronoi(\xb)}(\yb) = \yb$, and the Euclidean distance for that sample is zero.
Second, if the Voronoi predictor makes an incorrect prediction, then the term $\Vert \yb - \mathrm{proj}_{\backwardvoronoi(\xb)}(\yb) \Vert$ is non-zero, indicating the exact distance between the actual leak position and the predicted Voronoi cell.
The factor $1 / \text{(Number of Samples)}$ indicates that we are averaging over all samples, including both correct and incorrect predictions.
Alternatively, it may be of interest to know the error restricted to incorrect predictions only.
In that case, the only necessary change is replacing $1 / \text{(Number of Samples)}$ with $1/\text{(Number of Incorrect Predictions)}$.
In our evaluation below, we report both values.
The required projections are approximated by 100 iterations of Dykstra's projection algorithm \cite{boyle1986}.
The selected number of 100 iterations was tuned by hand and ensures convergence of the method according to our observations.

\subsection{Data cleaning}

\begin{figure}[t]
  \centering
  \includegraphics[trim=5 5 5 5, clip]{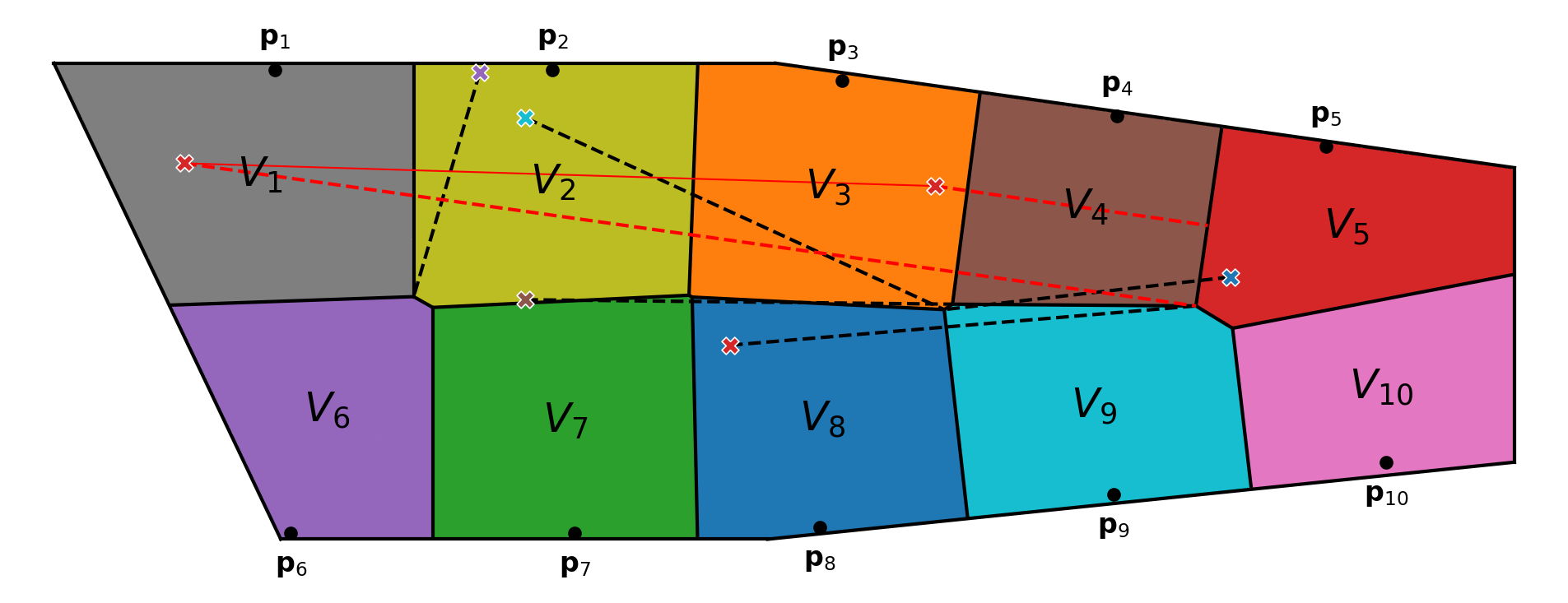}
  \caption{Results of outlier detection based on the Euclidean distance between the true leak positions and the classic Voronoi cell corresponding to the largest flow. Distances correspond to the lengths of the dashed lines (black for single-leak and red for multi-leak), and a sample is considered an outlier if the distance is greater than two meters.}
  \label{fig:outliers}
\end{figure}

There are several possible reasons why the data can contain outliers.
For example, a momentary malfunction of a flow meter, vacuum connection, or software can cause implausible measurement records.
Additionally, human error on the part of the experimenter can be a factor.
Collecting the data took approximately 250 hours, requiring high concentration.
Sources of human error include manually introduced leaks at positions different from the experimental design, improperly repaired leaks that cause biased measurements in subsequent steps, and typos or mix-ups when entering steady-state flows in a table.
Prior to evaluation, we aim to remove such outliers from the data.
To that end, we also use the Euclidean projection approach described above.
Specifically, if the distance between a true leak position and the associated predicted classic Voronoi cell is greater than a manually chosen threshold of two meters, then the respective sample is considered an outlier and removed from the data.
Figure~\ref{fig:outliers} shows the identified outliers, which include three single-leak samples and one two-leak sample.
For the two-leak samples, a sample is considered an outlier if both leaks are more than two meters from the corresponding classic Voronoi cell of the sensor with the largest flow.
For the sake of transparency, all evaluations were conducted on both the original and cleaned data, and we report all corresponding results below.

\subsection{Experimental results on single-leak data}

In this section, we present the experimental results obtained using classic and refined Voronoi predictors on single-leak data. As mentioned earlier, the original single-leak dataset contains 413 samples, five of which were identified as outliers during the cleaning process. Thus, the cleaned single-leak dataset contains 408 samples.
\medskip

\begin{figure}[t]
\centering

\begin{subfigure}{\linewidth}
    \centering
    \includegraphics[trim=5 5 5 5, clip]{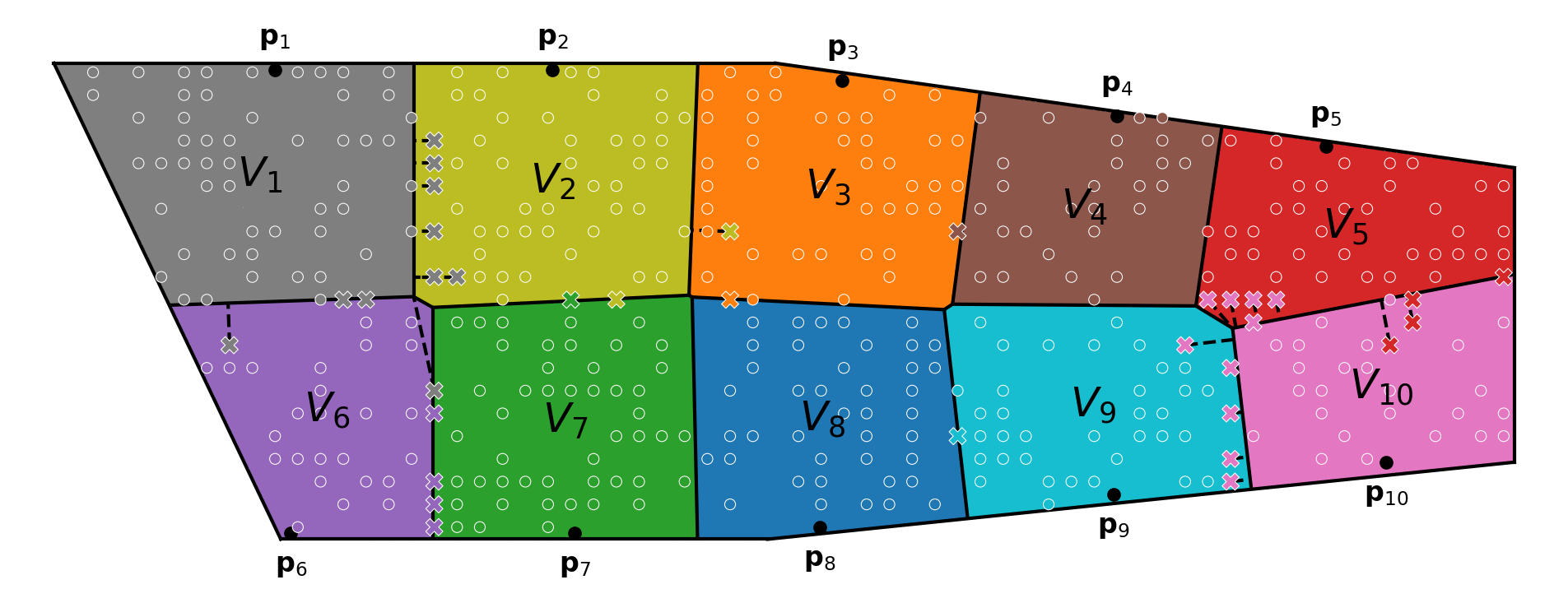}
    \caption{Classic Voronoi predictor.}
    \label{fig:evaluation_classic_single}
\end{subfigure}

\vspace{1em}

\begin{subfigure}{\linewidth}
    \centering
    \includegraphics[trim=5 5 5 5, clip]{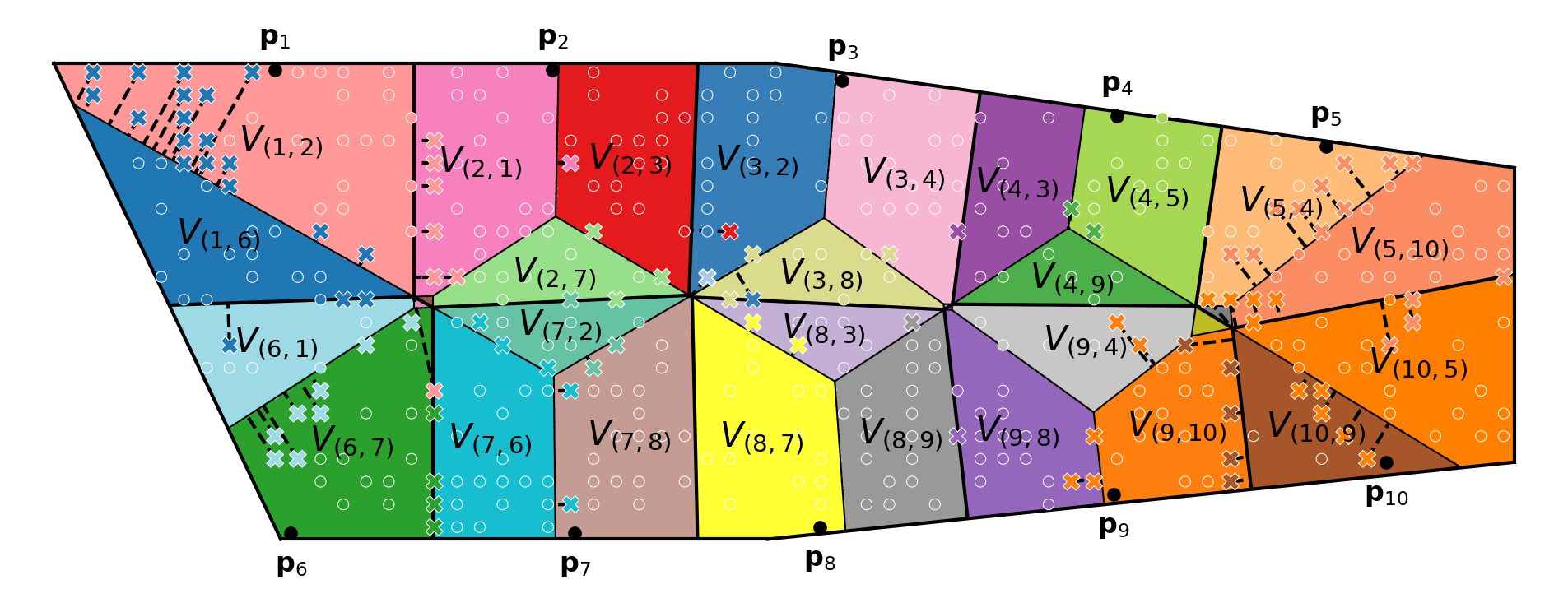}
    \caption{Refined Voronoi predictor.}
    \label{fig:evaluation_refined_single}
\end{subfigure}

\caption{Evaluation of Voronoi predictors on cleaned single-leak data. Each circle or cross indicates a leak position in the cleaned data. The circles represent correct predictions; that is, the actual leak position is within the Voronoi cell that was predicted based on the recorded flow data. Crosses represent incorrect predictions. Dashed lines connect the actual leak positions (crosses) and their associated projections on the predicted classic Voronoi cells.}
\label{fig:evaluation_single}
\end{figure}

Figure~\ref{fig:evaluation_single} provides an overview of the true leak locations in the cleaned single-leak dataset, along with the associated predicted Voronoi cells. Table~\ref{tab:experimental_results_single_leak} reports the corresponding accuracies and mean Euclidean distances, also for the original single-leak dataset. Clearly, according to the data cleaning procedure, both Voronoi predictors must perform better on the cleaned dataset in terms of both metrics. Here, we discuss the results on the cleaned single-leak dataset further.
\medskip

\begin{figure}[t]
    \centering

    \begin{subfigure}{.49\linewidth}
        \centering
        \includegraphics{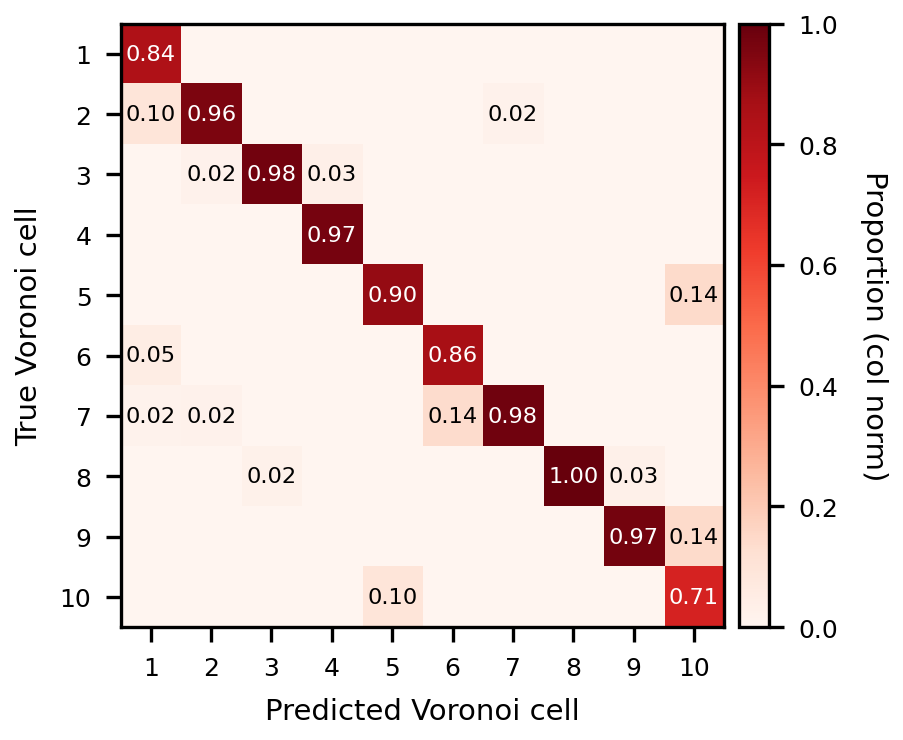}
        \caption{Column-wise normalization }
        \label{fig:confusion_matrix_classic_precision}
    \end{subfigure}\
    \begin{subfigure}{.49\linewidth}
        \centering
        \includegraphics{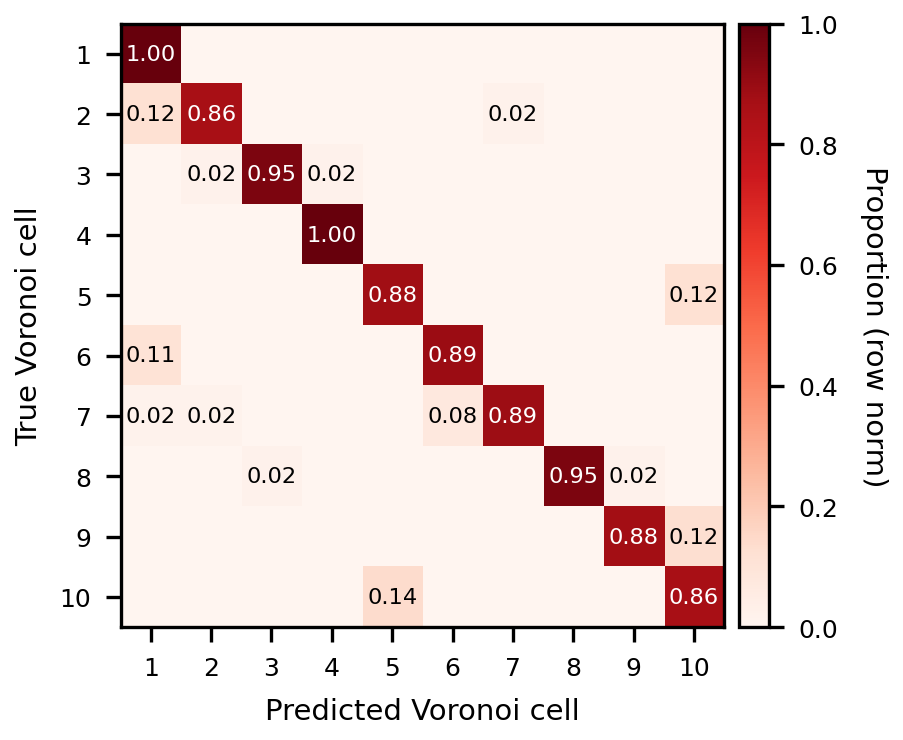}
        \caption{Row-wise normalization}
        \label{fig:confusion_matrix_classic_recall}
    \end{subfigure}
    
    \caption{Confusion matrices obtained with the classic Voronoi predictor on single-leak data, based on the frequency with which a predictor misclassifies pairs of classes in both directions, and normalized such that either columns (left) or rows (right) sum up to one.}
    \label{fig:confusion_matrices_classic}
\end{figure}

\begin{figure}
    \centering

    \begin{subfigure}{\linewidth}
        \centering
        \includegraphics{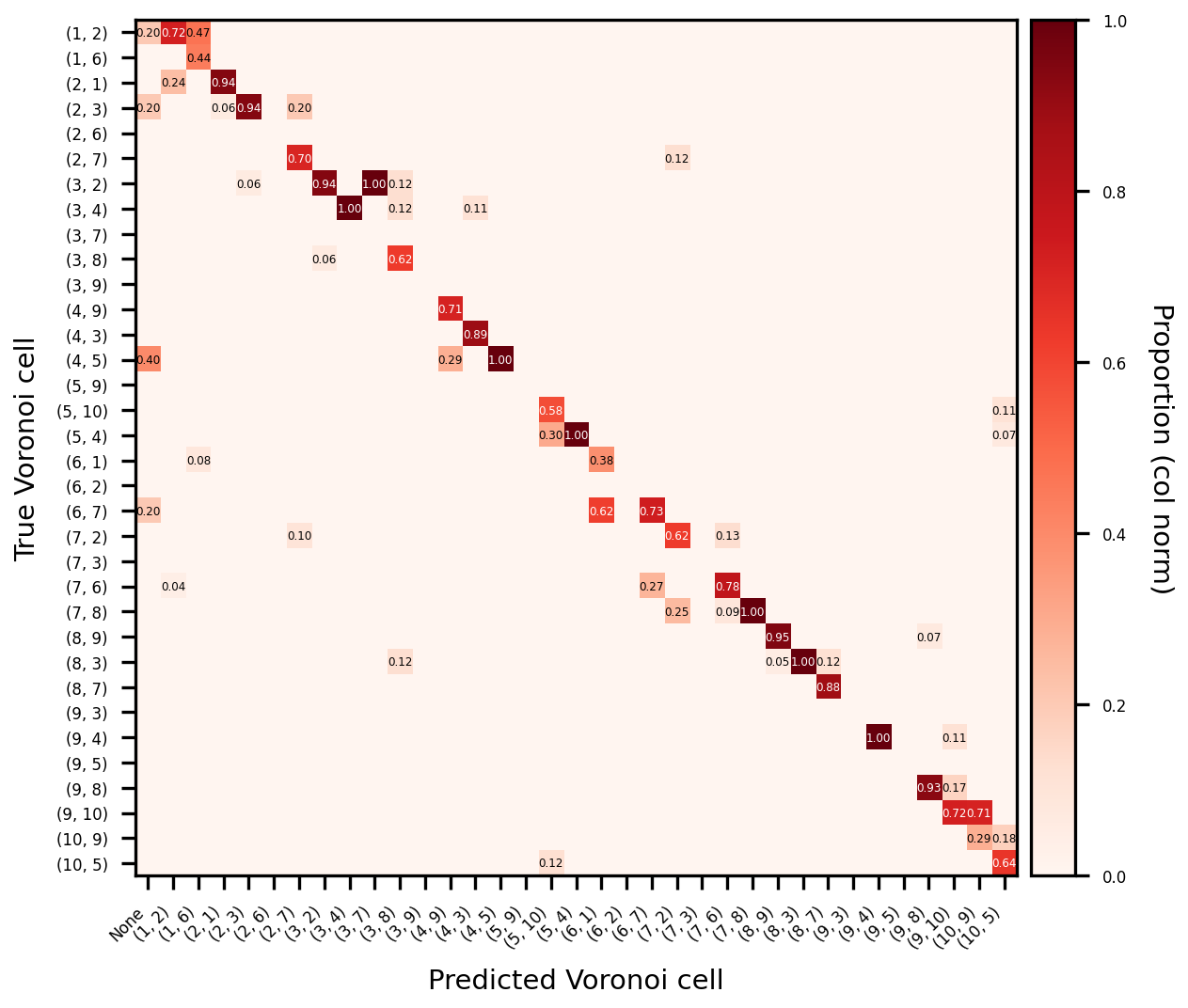}
        \caption{Row-wise normalization}
        \label{fig:confusion_matrix_refined_precision}
    \end{subfigure}\quad

    \vspace{1em}
    
    \begin{subfigure}{\linewidth}
        \centering
        \includegraphics{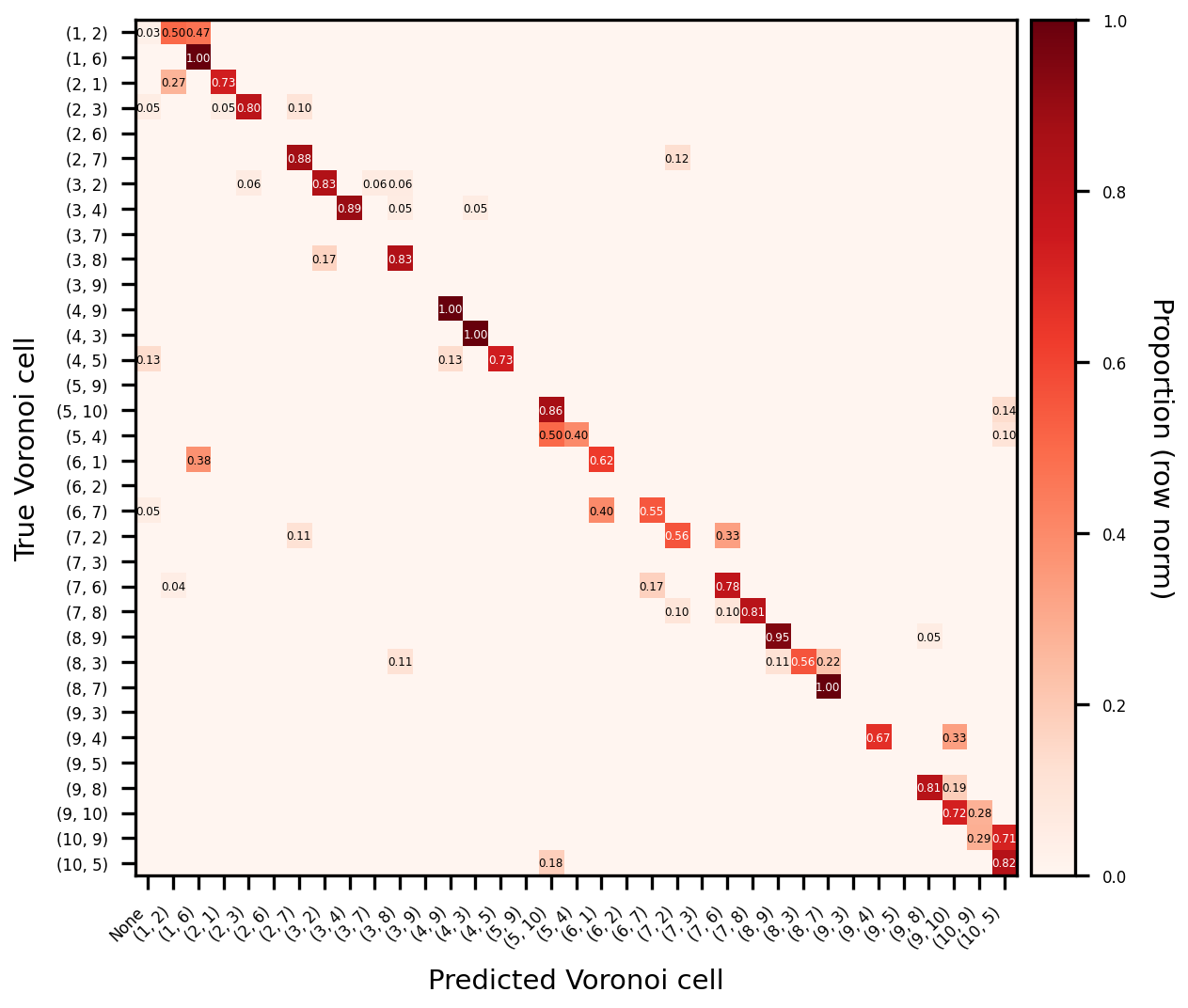}
        \caption{Column-wise normalization}
        \label{fig:confusion_matrix_refined_recall}
    \end{subfigure}
    
    \caption{Confusion matrices obtained with the refined Voronoi predictor on single-leak data.}
    \label{fig:confusion_matrices_refined}
\end{figure}

Figures~\ref{fig:confusion_matrices_classic} and \ref{fig:confusion_matrices_refined} illustrate \struc{confusion matrices} obtained using Voronoi predictors on cleaned single-leak data.
The entries in the matrices are based on the frequency with which a predictor misclassifies pairs of classes in both directions.
The absolute counts for each combination are then normalized either column-wise (Figures~\ref{fig:confusion_matrix_classic_precision} and \ref{fig:confusion_matrix_refined_precision}) or row-wise (Figures~\ref{fig:confusion_matrix_classic_recall} and \ref{fig:confusion_matrix_refined_recall}).
In the case of column normalization, the values on the diagonal of the confusion matrix correspond to class-specific precision, or the estimated probability that a prediction of a particular class is correct.
In the row-wise case, the diagonal values correspond to class-specific recall, or the estimated probability of a correct prediction given a particular true class.
\medskip

\begin{table}[b]
    \centering
    \caption{Evaluation on single-leak data}
    \label{tab:experimental_results_single_leak}
    \begin{tabular}{l l C{8.75em} C{8.75em} C{8.75em}}
        \toprule
        \makecell{Voronoi \\ type} &
        Dataset &
        \makecell{Accuracy \\ (\%)} &
        \makecell{Mean Euc. dist. \\ (full dataset, cm)} &
        \makecell{Mean Euc. dist. \\ (incorrect only, cm)} \\
        \midrule
        Classic & Original & 90.56 & 6.69 & 70.84 \\
        Classic & Cleaned & 91.67 & 1.71 & 20.58 \\
        \midrule
        Refined & Original & 73.85 & 11.53 & 45.69 \\
        Refined & Cleaned & 74.75 & 6.13 & 25.23 \\
        \bottomrule
    \end{tabular}
\end{table}

The classic Voronoi predictor achieves an accuracy of $91.67\ \%$. As Figure~\ref{fig:evaluation_classic_single} shows, most leak locations are marked with circles indicating correct predictions.
Most of the incorrectly predicted samples, marked with crosses, have leak locations that are only marginally outside the predicted classic Voronoi cell.
The mean Euclidean distance between the true leak locations and the predicted classic Voronoi cells for the incorrect predictions is $20.58$~cm.
This distance must be considered relative to the total part area, which is approximately $16 \times 5.2$~m.
False positive predictions concentrate on Voronoi cells $\voronoicell_1$ and $\voronoicell_{10}$, which lie in the upper left and lower right corners, respectively.
This can also be seen in the confusion matrix in Figure~\ref{fig:confusion_matrix_classic_precision}, where the first and last diagonal entries indicate low precision for these two cells.
Consequently, the classic Voronoi cells $\voronoicell_2$, $\voronoicell_5$, $\voronoicell_6$ and $\voronoicell_9$, which encompass the true leak locations corresponding to these false positive predictions, exhibit low recall, as indicated by the diagonal values in the confusion matrix in Figure~\ref{fig:confusion_matrix_classic_recall}.
The mean precision of all classic Voronoi cells is $91.66\ \%$ (the average of the diagonal entries in Figure~\ref{fig:confusion_matrix_classic_precision}).
The mean recall is $91.65\ \%$ (the average of the diagonal entries in Figure~\ref{fig:confusion_matrix_classic_recall}).
\medskip

The refined Voronoi predictor achieves an accuracy of $74.75\ \%$ on the single-leak data.
As with the classic Voronoi predictor, Figure~\ref{fig:evaluation_refined_single} shows that the largest deviations between the true leak locations and the predicted refined Voronoi cells occur at the left and right boundaries of the part's surface.
The confusion matrix in Figure~\ref{fig:confusion_matrix_refined_precision} shows that the precision values of these cells are significantly lower than the average of $72.94\ \%$.
The mean Euclidean distance with respect to incorrect predictions is $25.23$~cm.
Relative to the classic Voronoi predictor, this corresponds to a $22.59\ \%$ increase.
\medskip

\begin{figure}[t]
    \centering
    \includegraphics[trim=5 5 5 5, clip]{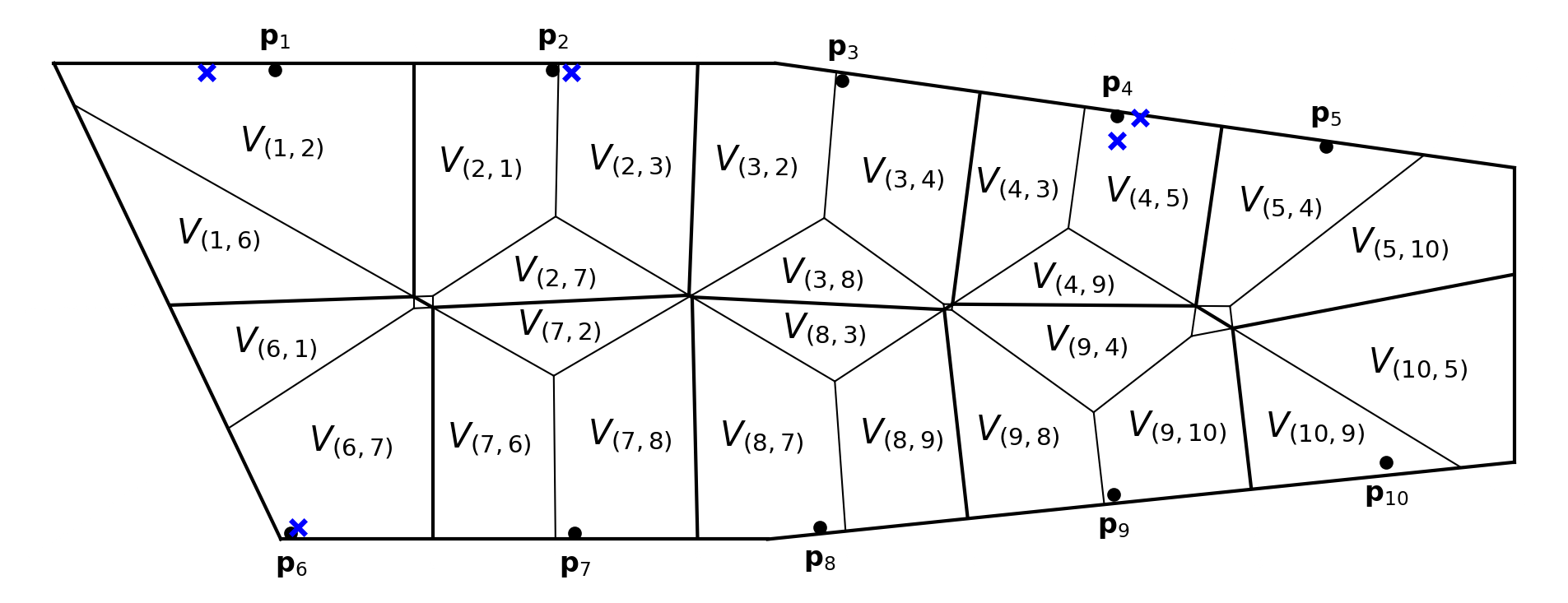}
    \caption{Invalid predictions of the refined Voronoi predictor. Each cross represents a case where the refined Voronoi predictor identifies a cell that does not exist in the given refined Voronoi diagram. Corresponding predictions are $\voronoicell_{(1,4)}$ (top left), $\voronoicell_{(2,10)}$ (top middle), $\voronoicell_{(4,10)}$ (top right, $2\times$), and $\voronoicell_{(6,10)}$ (bottom).}
    \label{fig:evaluation_refined_single_none}
\end{figure}

Furthermore, the confusion matrices in Figure~\ref{fig:confusion_matrices_refined} reveal a side effect of using refined Voronoi diagrams for leak detection:
On the one hand, there are Voronoi cells that are so small they do not contain any true leak positions and also are not predicted.
Specifically, these are $\voronoicell_{(2,6)}$, $\voronoicell_{(6,2)}$, $\voronoicell_{(3,7)}$, $\voronoicell_{(7,3)}$, $\voronoicell_{(3,9)}$, $\voronoicell_{(9,3)}$, $\voronoicell_{(5,9)}$, and $\voronoicell_{(9,5)}$, and the diagonal entries in the confusion matrices corresponding to these cells are zero.
Using $\voronoicell_{(2,6)}$ and $\voronoicell_{(6,2)}$ as an example, Figure~\ref{fig:evaluation_classic_single} shows that the corresponding classic Voronoi cells $\voronoicell_{2}$ and $\voronoicell_{6}$ are diagonally adjacent and have a relatively short common edge.
Consequently, the refined cells corresponding to their direct neighbors comprise most of the original area of $\voronoicell_{2}$ and $\voronoicell_{6}$.
On the other hand, in five cases, the refined Voronoi predictor identifies cells that do not exist in the given refined Voronoi diagram, namely $\voronoicell_{(2,10)}$, $\voronoicell_{(4,10)}$ (twice), $\voronoicell_{(6,10)}$, and $\voronoicell_{(1,4)}$.
These predictions, in the following referred to as \struc{invalid predictions}, are summarized in the \emph{None} category in Figure~\ref{fig:confusion_matrices_refined} (which has a corresponding column but no corresponding row in the confusion matrices because there are no corresponding true Voronoi cells) and illustrated in Figure~\ref{fig:evaluation_refined_single_none} (they are not contained in Figure~\ref{fig:evaluation_classic_single}). 
Notably, all corresponding true leak positions are located very close to one particular vacuum connection.

\subsection{Experimental results on multi-leak data}

\begin{table}[b]
\centering
\caption{Evaluation on multi-leak data}
\label{tab:experimental_results_multi_leak}
\begin{tabular}{l l
    C{4em} C{4em}
    C{4em} C{4em}
    C{4em} C{4em}}
\toprule
\makecell{Voronoi \\ type} &
Dataset &
\multicolumn{2}{c}{Accuracy (\%)} &
\multicolumn{2}{c}{\makecell{Mean Euc. dist. \\ (full dataset, cm)}} &
\multicolumn{2}{c}{\makecell{Mean Euc. dist. \\ (incorrect only, cm)}} \\
\cmidrule(lr){3-4}
\cmidrule(lr){5-6}
\cmidrule(lr){7-8}
& & Step 1 & Step 2 & Step 1 & Step 2 & Step 1 & Step 2 \\
\midrule
Classic & Original & 92.95 & 87.04 & 1.97 & 4.42 & 27.93 & 34.12 \\
Classic & Cleaned  & 93.06 & 88.31 & 1.30 & 1.30 & 18.76 & 11.14 \\
\midrule
Refined & Original & 33.86 & 70.47 & 28.89 & 12.54 & 65.07 & 43.15 \\
Refined & Cleaned  & 33.80 & 71.23 & 29.20 & 6.07 & 65.63 & 21.45 \\
\bottomrule
\end{tabular}
\end{table}

\begin{figure}[t]
\centering

\begin{subfigure}{\linewidth}
    \centering
    \includegraphics[trim=5 5 5 5, clip]{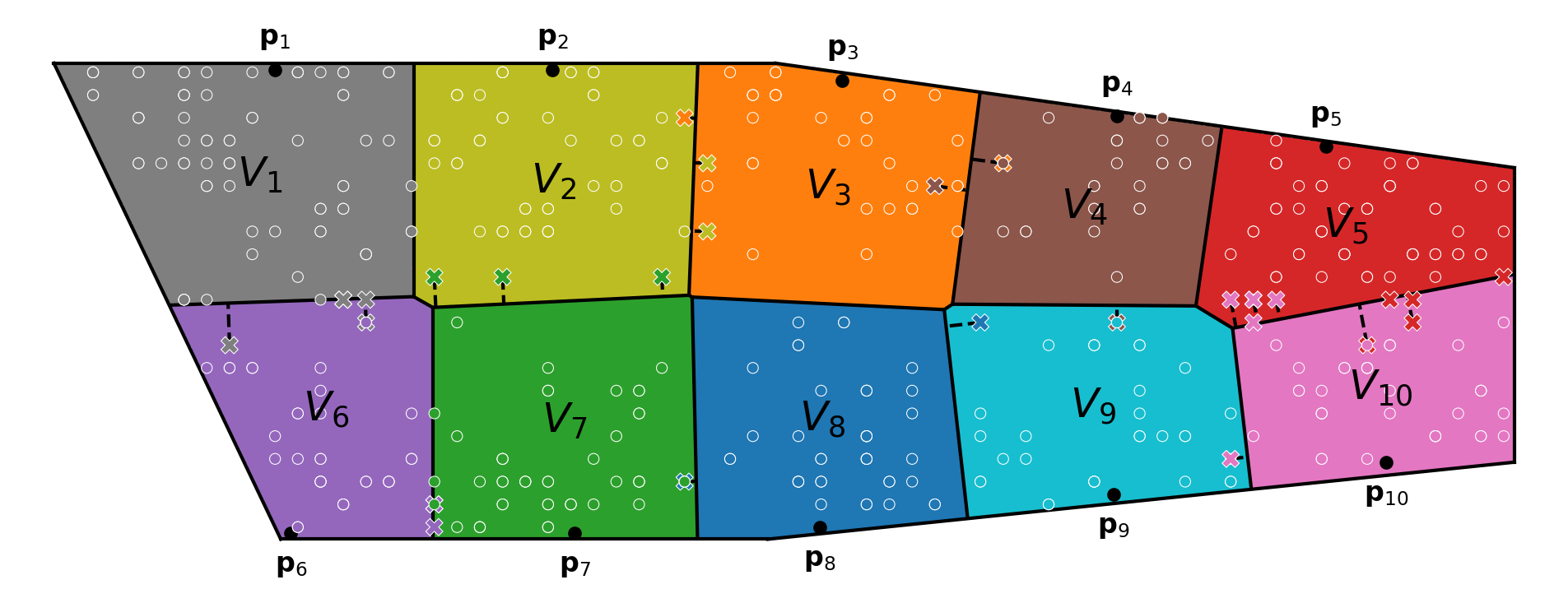}
    \caption{Classic Voronoi predictor.}
    \label{fig:evaluation_classic_multi}
\end{subfigure}

\vspace{1em}

\begin{subfigure}{\linewidth}
    \centering
    \includegraphics[trim=5 5 5 5, clip]{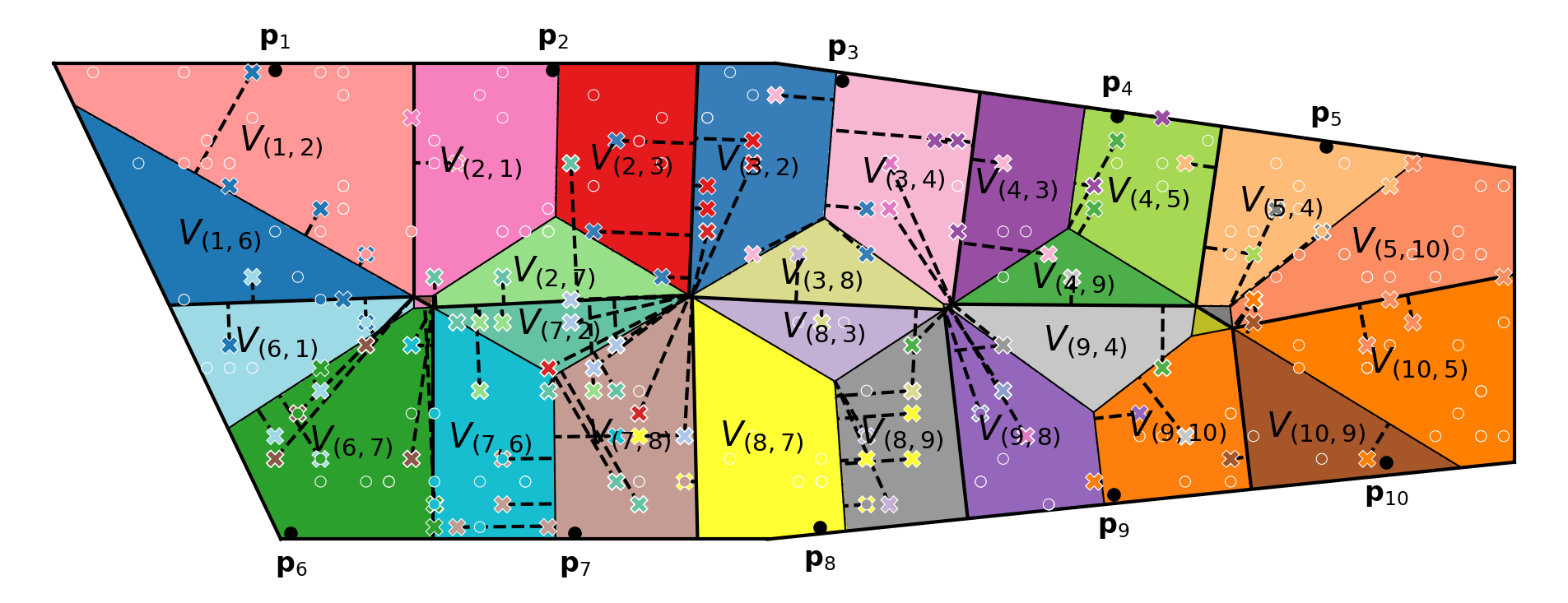}
    \caption{Refined Voronoi predictor.}
    \label{fig:evaluation_refined_multi}
\end{subfigure}

\caption{Evaluation of Voronoi predictors on cleaned multi-leak data. Each circle or cross indicates a leak position in the cleaned data. The circles represent a correct prediction in step~1 of the repeated leak localization approach, meaning at least one of the two true leak positions in the same sample is within the predicted cell. The crosses represent cases in which neither true leak position is within the predicted cell in step one. Only the leak position with the smaller distance to the predicted cell and its projection onto the predicted cell are shown, connected by a dashed line.}
\label{fig:evaluation_multi}
\end{figure}

Now, we evaluate the use of classic and refined Voronoi predictors with multi-leak data. By definition, both predictors are designed to detect single leaks.
In Section~\ref{sec:strategies_for_multi_leak_localization}, we introduced the repeated leak localization strategy to expand the capabilities of our predictors for multi-leak detection.
The original dataset contains 440 two-leak examples. One example was removed during data cleaning of the two-leak data.
Seven more examples were removed because they included leak positions that were removed during the cleaning of the single-leak data.
This left 432 examples in the cleaned two-leak dataset. After these steps, each position of a leak in a two-leak example also appears as a single-leak example.
This enables us to simulate the repeated leak detection strategy as follows:
First, we apply a Voronoi predictor to the flow values induced by the two leaks and evaluate the predictor's ability to identify one of the two leaks.
If successful, we feed the flow induced by the remaining leak to the same Voronoi predictor and evaluate the second step. 
\medskip

Table~\ref{tab:experimental_results_multi_leak} summarizes the results for multi-leak data.
Step~1 accuracy is defined as the percentage of samples for which at least one true leak position falls within the predicted cell.
Step~2 accuracy is defined as the percentage of samples for which the second true leak position falls within the predicted cell.
The step~1 Euclidean distance is the minimum distance between the two leak positions and the predicted cell, i.e., it is non-zero only if both true leak positions are outside the predicted cell.
Both step~2 accuracy and Euclidean distance are computed only among samples for which step~1 is successful.
Figure~\ref{fig:evaluation_multi} illustrates the step~1 predictions for both Voronoi predictors.
Again, the crosses indicate incorrect predictions, i.e., that none of the two true leak positions in a sample is within the predicted cell in step one.
For each incorrectly predicted sample, only the leak position with the smaller distance to the predicted cell is shown.
\medskip

The classic Voronoi predictor achieves $93.06\ \%$ accuracy in step~1, even surpassing the accuracy on single-leak data.
Step~2 accuracy is $88.31\ \%$.
Likewise, the mean Euclidean distances are even lower in the two-leak case than in the single-leak case.
At first glance, Figures~\ref{fig:evaluation_classic_single} and \ref{fig:evaluation_classic_multi} appear similar.
However, upon closer inspection, it becomes apparent that the density of the marked leak positions in Figure~\ref{fig:evaluation_classic_multi} is significantly lower toward the horizontal center of the part.
The missing leak positions correspond to incorrect step~1 predictions and are those with the greater distance to the predicted cell among the two leak positions in their sample.
\medskip

The refined Voronoi predictor achieves an accuracy of $33.80\ \%$ in step~1 and $71.23\ \%$ in step~2.
The significant discrepancy between the accuracy of the refined Voronoi predictor in the single-leak case and the accuracy in step~1 of the two-leak case is primarily due to invalid predictions, i.e., predictions for which the predicted refined Voronoi cell does not exist in the diagram.
Of the 432 examples in the cleaned two-leak data, the refined Voronoi predictor produces 146 correct predictions in step~1 (corresponding to the $33.80\ \%$ accuracy), 169 invalid predictions ($39.12\ \%$, not included in Figure~\ref{fig:evaluation_refined_multi}), and 117 valid but incorrect predictions ($27.08\ \%$).
Since invalid predictions do not correspond to any cell in the refined Voronoi diagram, the respective samples were excluded from the calculation of the reported mean Euclidean distances.
Among the valid predictions, the mean Euclidean distances of the refined Voronoi predictor in step~1 are significantly larger than the mean Euclidean distances in the single-leak case.
This can also be seen by comparing Figures~\ref{fig:evaluation_refined_single} and \ref{fig:evaluation_refined_multi}.
\medskip

Before concluding the section on our experiments, we investigate the observed invalid predictions further.
We hypothesize that these predictions may contain useful information about the spatial distribution of multiple leaks, which the refined Voronoi predictor is currently unable to exploit.
This is because the predictor establishes a one-to-one correspondence between the indices of the two vacuum connections with the largest measurements and the corresponding refined Voronoi cell.
If a corresponding refined Voronoi cell does not exist, the prediction cannot be used to draw conclusions about leak locations.
One possible solution is to relax the aforementioned one-to-one correspondence where appropriate.
While developing further prediction models is beyond the scope of this paper, the following analysis aims to provide initial insights and directions for future development.
\medskip

\begin{figure}
    \centering
    \includegraphics[trim=5 5 5 5, clip]{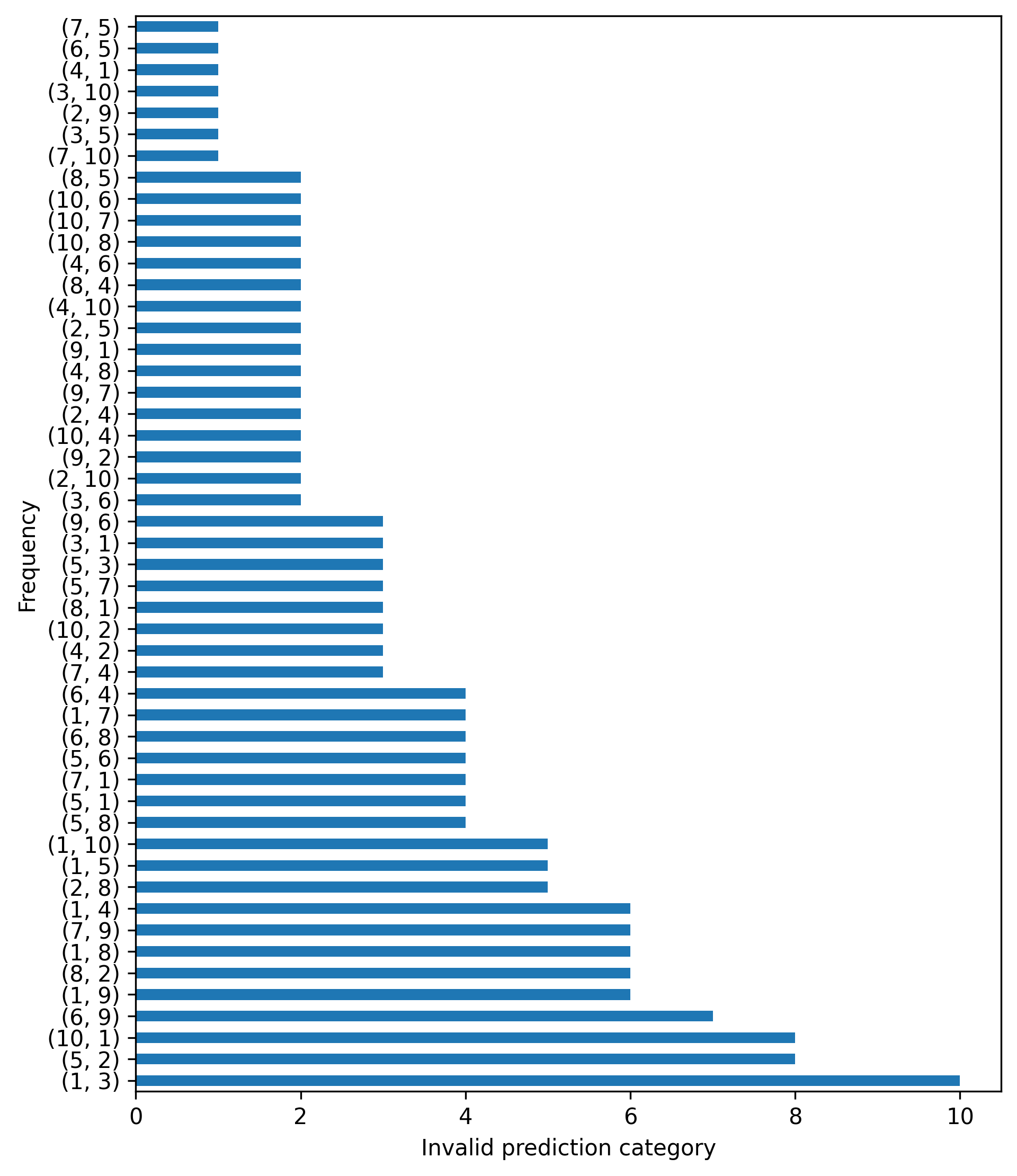}
    \caption{Histogram of invalid prediction categories for the refined Voronoi predictor.}
    \label{fig:invalid_predictions_histogram}
\end{figure}

Figure~\ref{fig:invalid_predictions_histogram} shows a histogram of the 50 invalid prediction categories that occur when the refined Voronoi predictor is applied to the cleaned two-leak data.
The related Voronoi diagram has 34 cells, which corresponds to the number of valid prediction categories.
Since there are 90 possible orderings of two out of ten measurements, only six theoretically possible invalid prediction categories do not occur.
These missing categories are $(4, 7)$, $(6, 3)$, $(6, 10)$, $(8, 6)$, $(8, 10)$, and $(10, 3)$.
\medskip

\begin{figure}[t]
\centering

\begin{subfigure}{\linewidth}
    \centering
    \includegraphics[trim=5 5 5 5, clip]{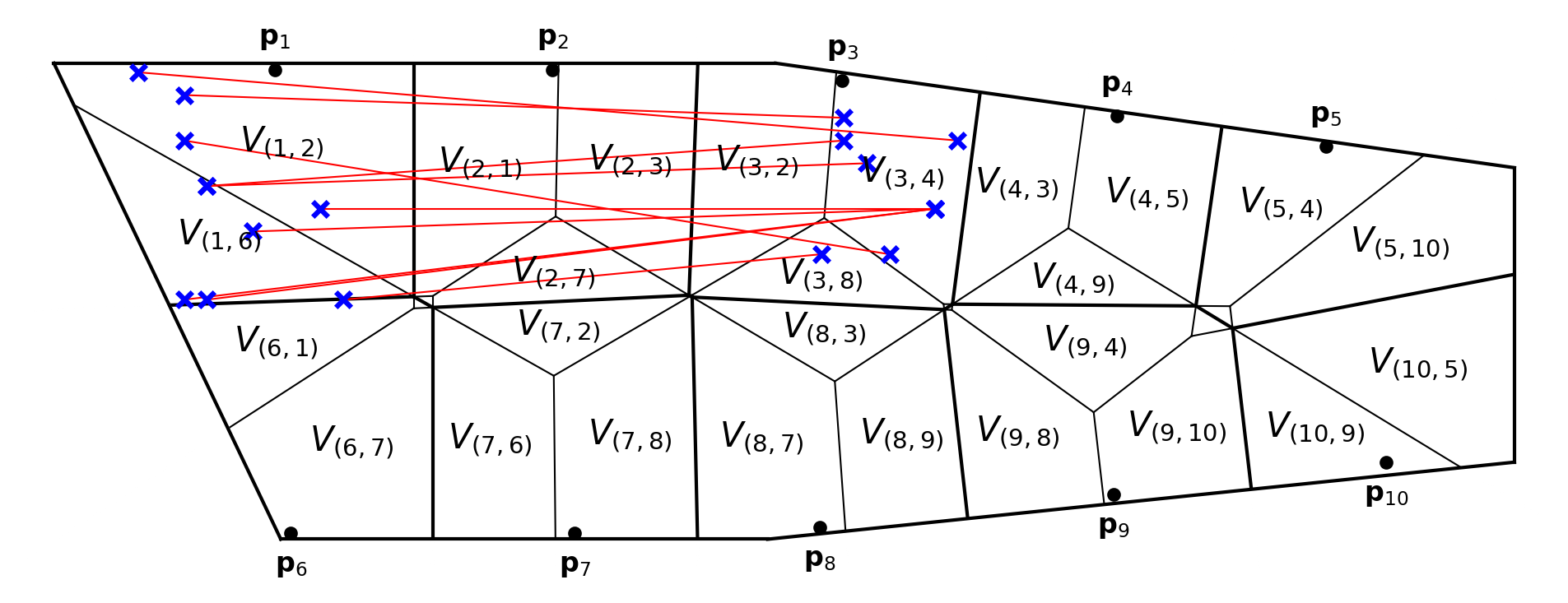}
    \caption{Invalid prediction $(1, 3)$}
    \label{fig:invalid_predictions_1_3}
\end{subfigure}

\vspace{1em}

\begin{subfigure}{\linewidth}
    \centering
    \includegraphics[trim=5 5 5 5, clip]{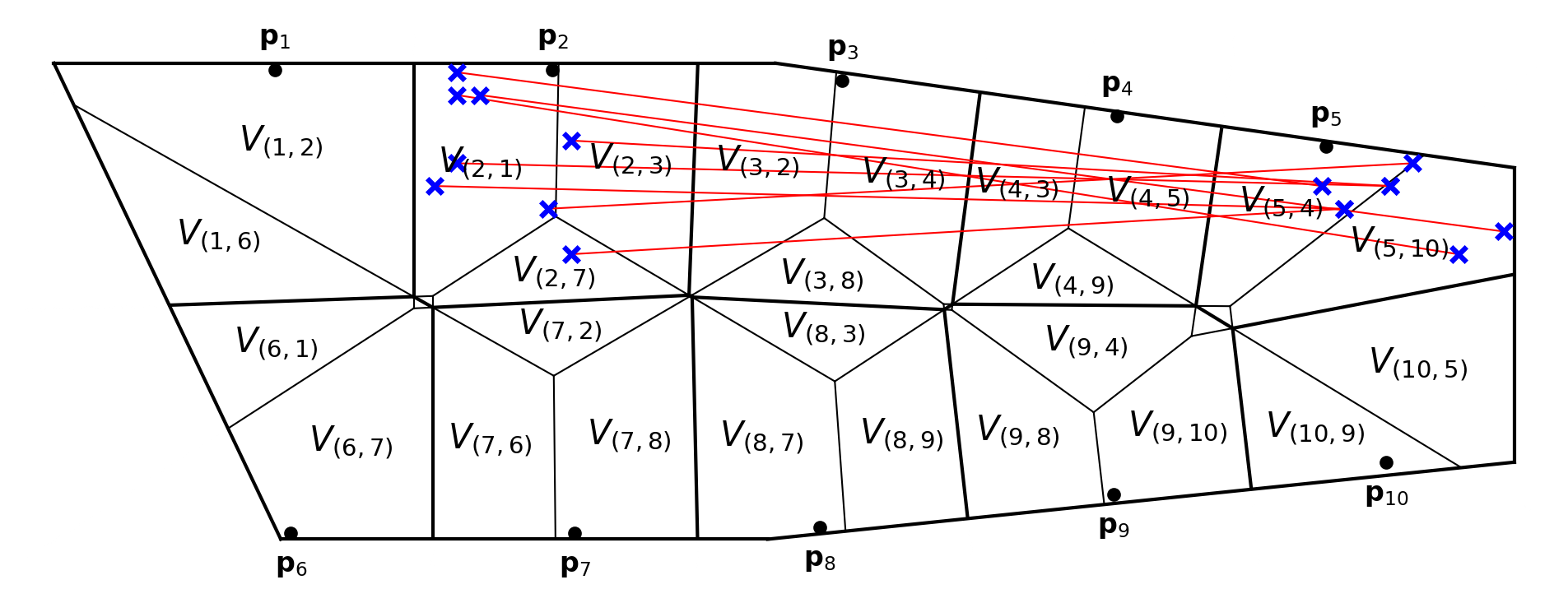}
    \caption{Invalid prediction $(5, 2)$}
    \label{fig:invalid_predictions_5_2}
\end{subfigure}

\vspace{1em}

\begin{subfigure}{\linewidth}
    \centering
    \includegraphics[trim=5 5 5 5, clip]{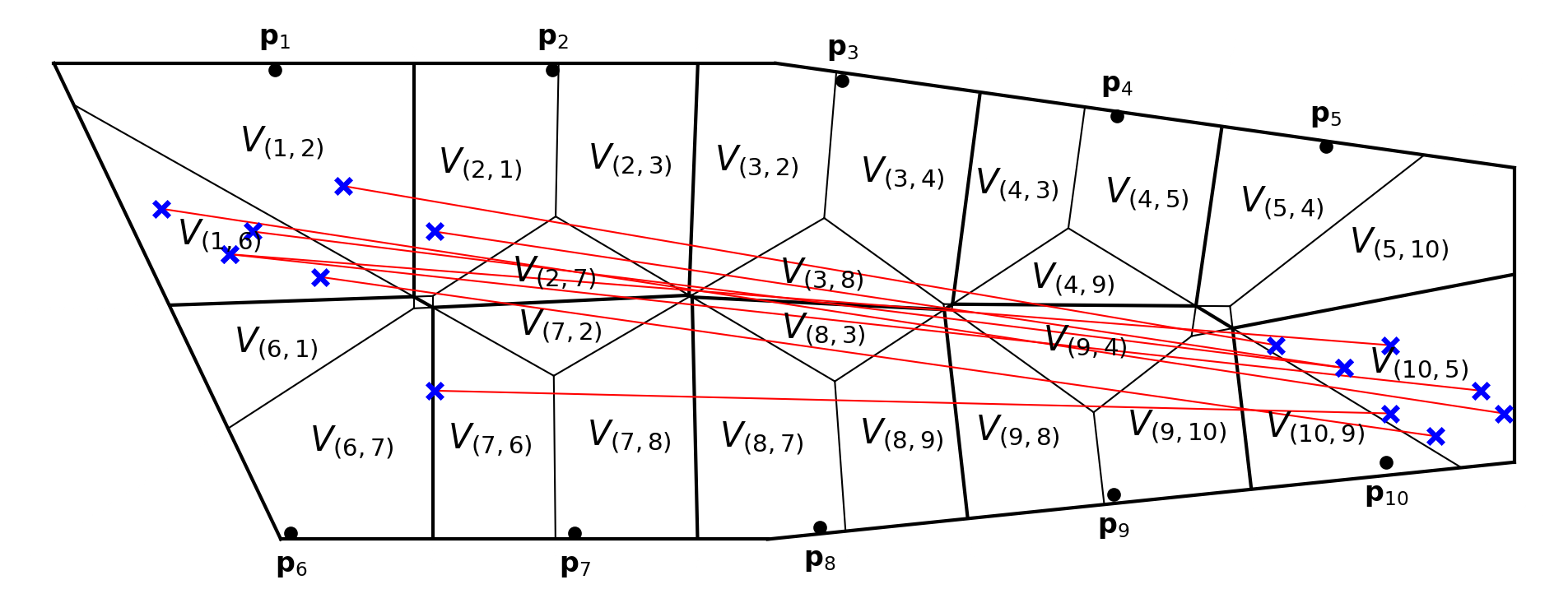}
    \caption{Invalid prediction $(10, 1)$}
    \label{fig:invalid_predictions_10_1}
\end{subfigure}

\caption{Spatial distribution of leak locations for the three most frequent invalid prediction categories of the refined Voronoi predictor.}
\label{fig:invalid_predictions}
\end{figure}

Figure~\ref{fig:invalid_predictions} shows the spatial distribution of leak locations associated with the three most frequent invalid prediction categories $(1, 3)$, $(5, 2)$, and $(10, 1)$.
The leak locations in Figure~\ref{fig:invalid_predictions_1_3} show that all two-leak samples associated with the prediction $(1, 3)$ feature one leak location in the classic Voronoi cell $\voronoicell_1$ (except in one case where the leak location is in $\voronoicell_6$ but very close to the shared edge of $\voronoicell_1$ and $\voronoicell_6$) and one leak location in the classic Voronoi cell $\voronoicell_3$.
Moreover, it can be observed that most leak locations even concentrate in or close to refined Voronoi cells $\voronoicell_{(1, 6)}$ and $\voronoicell_{(3, 4)}$.
Note that two of the leak locations in $\voronoicell_3$ appear in two samples which is why there are nine locations in $\voronoicell_1$ but only seven locations in $\voronoicell_3$.
Basically, the same statements hold for the leak locations illustrated in Figures~\ref{fig:invalid_predictions_5_2} and \ref{fig:invalid_predictions_10_1}.

\section{Discussion}
\label{sec:discussion}

Here, we discuss the experimental results reported in Section~\ref{sec:experiments} and their potential implications for practical leak detection.
For the sake of brevity, we focus on the main aspects.

\subsubsection*{Classic Voronoi predictor}

Our results indicate that the classic Voronoi predictor reliably performs on both one-leak and two-leak data.
In practice, the restriction of leak positions to classic Voronoi cells already facilitates subsequent techniques, such as ultrasonic, helium, or thermographic leak tests a lot.
In the event of an incorrect prediction, the true leak position is usually very close to the predicted cell (on average, less than 20 cm).
Therefore, most incorrect predictions can be accounted for by introducing a small safety margin to the predicted classic Voronoi cell.
    
\subsubsection*{Refined Voronoi predictor}
The results obtained with the refined Voronoi predictor require further differentiation.
This is particularly because of the significant discrepancy between the results for one-leak and two-leak data.

\begin{description}
    \item[One-leak data] The accuracy is lower than that of the classic Voronoi predictor, which is expected because of the smaller cells.
    However, the same statement as above holds true:
    A small security margin around the predicted refined Voronoi cells captures a significant portion of all incorrect predictions as well (the mean distance is around 25 cm).
    For practical applications, this indicates that a prioritized search within and around a refined Voronoi cell is effective before extending the search area to a classic Voronoi cell if necessary.
    \item[Two-leak data] While the classic Voronoi predictor has similar accuracy in steps one and two as in the single-leak case, this does not hold for the refined Voronoi predictor.
    The refined Voronoi predictor's step~1 accuracy is only $33.8\ \%$, which is less than half of its step 2 accuracy.
    This suggests that the refined Voronoi predictor has limited practical use with the repeated prediction strategy.
\end{description}

\subsubsection*{Invalid predictions}

The poor accuracy of the refined Voronoi predictor in step~1 of the repeated strategy is largely due to invalid predictions, which account for 169 out of 432 samples ($39.12 \%$; the step~1 accuracy among only valid predictions is $55.52\ \%$).
However, our investigation of the most frequently observed invalid prediction categories suggests that these predictions contain useful information that the refined Voronoi predictor just cannot exploit.
First, our results indicate that invalid predictions are a relatively robust indicator of the presence of multiple leaks.
Specifically, only five invalid predictions were observed in the one-leak case, versus 169 invalid predictions in the two-leak case.
Furthermore, in the two-leak case, it is often even possible to correctly predict the classic Voronoi cells containing the leaks in a single step.
For 144 of the 169 invalid predictions ($85.21\ \%$), the predicted tuple $(i_1, i_2)$ points to classic Voronoi cells $\voronoicell_{i_1}$ and $\voronoicell_{i_2}$ which both contain a leak.
For additional 20 of the 169 invalid predictions ($11.83\ \%$), only $\voronoicell_{i_1}$ contains a leak.
And in 5 of 169 cases ($2.96\ \%$), only $\voronoicell_{i_2}$ contains a leak.
The associated search strategy is the simultaneous search strategy briefly discussed above.
Unlike the repeated search strategy, the simultaneous strategy only requires a single measurement.
The simultaneous strategy can also be applied beyond invalid predictions of the refined Voronoi predictor.
Overall, it locates both leaks correctly for 311 of 432 two-leak samples ($72.00\ \%$).
For 120 of the 432 two-leak samples ($27.78\ \%$), it correctly locates at least one leak.
In only one case out of 432, neither of the cells $\voronoicell_{i_1}$ and $\voronoicell_{i_2}$ contains a leak.
However, as also mentioned above, the use of the simultaneous strategy is limited when the number of leaks is unknown.
One takeaway is that if step~1 of the repeated strategy is unsuccessful, it makes sense to look at the cell corresponding to the second-largest measurement.
\bigskip

We summarize our discussion as follows:
\begin{itemize}[itemsep=6pt]
    \item The classic Voronoi predictor combined with the repeated leak localization strategy reliably predicts a cell containing a leak in the first shot with an accuracy of approximately $92\ \%$.
    This reveals an effective strategy for both single and multiple leaks localization reducing the search area by a factor of about ten (number of vacuum connections).
    \item If the classic Voronoi cell associated with the largest flow measurement does not contain a leak, then the cell corresponding to the second largest flow measurement likely does.
    \item The standalone refined Voronoi predictor is less reliable but useful in a prioritized search strategy combined with the classic Voronoi predictor.
    \item Invalid predictions are a reliable indicator of multiple leaks.
    \item We expect the underlying information about the order of measurements to be a useful feature for follow-up localization strategies, e.g., based on neural networks.
\end{itemize}

\section{Conclusion}
\label{sec:conclusion}

In this paper, we investigated the use of Voronoi diagrams for leakage detection in vacuum setups for composite manufacturing.
In addition to the classic Voronoi diagram, which is based on the positions of individual vacuum connections, we introduced a refined version that considers the order of distances from a point to multiple connections.
We derived knowledge-based predictors from both variants that can detect single or multiple leaks with a suitable repeated search strategy.
These predictors were evaluated and compared using one- and two-leak scenarios generated on an industrial-scale wing cover mold.
The classic Voronoi predictor reliably narrows down the search area for a leak from the entire mold surface to a single Voronoi cell.
While the refined Voronoi predictor is less precise as a standalone predictor, our analysis revealed that the underlying order of flow measurements contains important information that can enable multi-leak detection and enhance predictive models in future work.
One potential enhancement is the use of data-based strategies. First, such strategies can adapt Voronoi diagrams to experimental data, i.e., shift nodes and edges slightly to transform incorrect predictions into correct ones.
Second, classification models based on flow data and Voronoi diagrams can likely improve upon the knowledge-driven predictors presented in this work, especially when the number of leaks is unknown.
Finally, an important industrial use case and future direction is leakage detection on surfaces that incorporate flow barriers or additional flow aids below the vacuum bag.
Both prevent isotropic flow.
Examples of this include setups with stringers or other reinforcing elements, additional hoses, and areas without breather material.
We consider this work to be the foundation for the aforementioned future directions that are the subject of current activities or will be part of future research.

\clearpage
\bibliographystyle{amsalpha}
\bibliography{./leakage_detection_voronoi_arxiv_v1}

@inproceedings{shamos1975closest,
  author       = {Shamos, Michael Ian and Hoey, Dan},
  title        = {Closest-Point Problems},
  booktitle    = {Proceedings of the 16th Annual Symposium on Foundations of Computer Science (FOCS)},
  pages        = {151--162},
  organization = {IEEE},
  year         = {1975}
}

@book{aurenhammer2013voronoi,
  author    = {Aurenhammer, Franz and Klein, Rolf and Lee, Der-Tsai},
  title     = {Voronoi Diagrams and {D}elaunay Triangulations},
  publisher = {World Scientific},
  year      = {2013}
}

@inproceedings{boyle1986,
  author       = {Boyle, James P. and Dykstra, Richard L.},
  title        = {A Method for Finding Projections onto the Intersection of Convex Sets in {H}ilbert Spaces},
  booktitle    = {Advances in Order Restricted Statistical Inference},
  pages        = {28--47},
  organization = {Springer},
  year         = {1986}
}

@inproceedings{dlrk2011,
  author    = {Ucan, Hakan and B{\"o}lke, Jens and Krombholz, Christian and Gobbi, Henrique and Meyer, Mathias},
  title     = {Robotergest{\"u}tzte {L}eckageerkennung an {V}akuumaufbauten mittels {T}hermografie},
  booktitle = {Deutscher Luft- und Raumfahrtkongress},
  year      = {2011}
}

@book{goodfellow2016,
  author    = {Goodfellow, Ian and Bengio, Yoshua and Courville, Aaron},
  title     = {Deep Learning},
  publisher = {MIT Press},
  year      = {2016}
}

@article{haschenburger2019a,
  author    = {Haschenburger, Anja and Heim, Clemens},
  title     = {Two-Stage Leak Detection in Vacuum Bags for the Production of Fibre-Reinforced Composite Components},
  journal   = {CEAS Aeronautical Journal},
  volume    = {10},
  number    = {3},
  pages     = {885--892},
  publisher = {Springer},
  year      = {2019}
}

@article{haschenburger2021a,
  author    = {Haschenburger, Anja and Menke, Niklas and St{\"u}ve, Jan},
  title     = {Sensor-Based Leakage Detection in Vacuum Bagging},
  journal   = {The International Journal of Advanced Manufacturing Technology},
  volume    = {116},
  number    = {7},
  pages     = {2413--2424},
  publisher = {Springer},
  year      = {2021}
}

@inproceedings{brauer2022a,
  author       = {Brauer, Christoph and Lorenz, Dirk and Tondji, Lionel},
  title        = {Group Equivariant Networks for Leakage Detection in Vacuum Bagging},
  booktitle    = {European Signal Processing Conference (EUSIPCO)},
  pages        = {1437--1441},
  organization = {IEEE},
  year         = {2022}
}

@article{haschenburger2022computational,
  author    = {Haschenburger, Anja and Onorato, L. and Sujahudeen, M. S. and Taraczky, D. S. and Osis, A. and Bracke, A. R. S. and Byelov, M. D. and Vermeulen, F. I. and Oosthoek, E. H. Q.},
  title     = {Computational Methods for Leakage Localisation in a Vacuum Bag Using Volumetric Flow Rate Measurements},
  journal   = {Production Engineering},
  volume    = {16},
  number    = {6},
  pages     = {823--835},
  publisher = {Springer},
  year      = {2022}
}

@phdthesis{dlr190099,
  author = {Haschenburger, Anja},
  title  = {Influence and Detection of Vacuum Bag Leakages in Composites Manufacturing},
  school = {Delft University of Technology},
  month  = {November},
  year   = {2022},
  doi    = {10.4233/uuid:03637286-f682-4602-9890-2ac1c0102599},
  url    = {https://elib.dlr.de/190099/}
}

@book{Ziegler,
  author    = {Ziegler, G{\"u}nter M.},
  title     = {Lectures on Polytopes},
  publisher = {Springer},
  address   = {New York},
  year      = {2007}
}

@book{Gruenbaum,
  author    = {Gr{\"u}nbaum, Branko and Klee, Victor and Perles, Micha A. and Shephard, Geoffrey C.},
  title     = {Convex Polytopes},
  publisher = {Springer},
  volume    = {16},
  year      = {1967}
}

@mastersthesis{naveenachandran2023,
  author = {Naveenachandran, Sreerag V.},
  title  = {Data-Based Leakage Detection and Uncertainty Quantification in the Manufacturing of Large-Scale CFRP Components},
  school = {TU Braunschweig},
  year   = {2023}
}

@inproceedings{ChazelleEdelsbrunner,
  author    = {Chazelle, Bernard and Edelsbrunner, Herbert},
  title     = {An Improved Algorithm for Constructing k-th Order {V}oronoi Diagrams},
  booktitle = {Proceedings of the First Annual Symposium on Computational Geometry},
  pages     = {228--234},
  year      = {1985}
}

@inproceedings{Chan2024OptimalHOVD,
  author       = {Chan, Timothy M. and Cheng, Pingan and Zheng, Da Wei},
  title        = {An Optimal Algorithm for Higher-Order {V}oronoi Diagrams in the Plane: The Usefulness of Nondeterminism},
  booktitle    = {Proceedings of the Annual ACM-SIAM Symposium on Discrete Algorithms (SODA)},
  pages        = {4451--4463},
  organization = {SIAM},
  year         = {2024}
}

@article{Indermitte:VoronoiBio,
  author    = {Indermitte, Claude and Liebling, Thomas M. and Troyanov, Marc and Cl{\'e}men{\c{c}}on, Heinz},
  title     = {Voronoi Diagrams on Piecewise Flat Surfaces and an Application to Biological Growth},
  journal   = {Theoretical Computer Science},
  volume    = {263},
  number    = {1--2},
  pages     = {263--274},
  publisher = {Elsevier},
  year      = {2001}
}

@article{Brassel:ThiessenPolygons,
  author    = {Brassel, Kurt E. and Reif, Douglas},
  title     = {A Procedure to Generate Thiessen Polygons},
  journal   = {Geographical Analysis},
  volume    = {11},
  number    = {3},
  pages     = {289--303},
  publisher = {Wiley},
  year      = {1979}
}

@article{Blatov:VoronoiCrystal,
  author    = {Blatov, Vladislav A.},
  title     = {Voronoi--{D}irichlet Polyhedra in Crystal Chemistry: Theory and Applications},
  journal   = {Crystallography Reviews},
  volume    = {10},
  number    = {4},
  pages     = {249--318},
  publisher = {Taylor \& Francis},
  year      = {2004}
}

\clearpage
\appendix

\section{Discrete Geometry}
\label{appendix:discrete_geometry}

We review some basic concepts and notations from discrete geometry. 
Fur further details, see introductory sources such as \cite{Gruenbaum,Ziegler}.
\medskip

Let $A$ be an arbitrary subset of $\R^n$.
The \struc{convex hull} of $A$ is given as

\begin{align*}
	\struc{\conv(A)} \ := \ \left\{ \sum_{i=1}^{k}{\lambda_i \ab_i} \ \middle| \ \ab_i \in A, \ \lambda_i \geq 0, \ \sum_{i=1}^{k}{\lambda_i}=1, \ k\in\N \right\}.
\end{align*}
\smallskip

\noindent Moreover, $A$ is called \struc{convex} if $\conv(A) = A$.
Geometrically, $A$ is convex if and only if it holds for every line segment connecting two points $\Vector{a},\Vector{b} \in A$ that the entire line segment $[\Vector{a},\Vector{b}]$ is also contained in $A$.
\medskip

A dual approach for constructing convex sets is given by the intersection of half spaces.
An $(n-1)$-dimensional subspace $H$ of $\R^n$ is called an \struc{(affine) hyperplane}.
$H$ is uniquely determined by a linear form, i.e., 

\begin{align*}
	H \ = \ \left\{\xb \in \R^n \suchthat \langle \Vector{a},\xb \rangle = b\right\}
\end{align*}
\smallskip

\noindent for some $\Vector{a} \in \R^n \setminus \{\Vector{0}\}$ and $b\in \R$.
The vector $\Vector{a}\in \R^n$ in this representation is unique up to nonzero multiples. 
In fact, $\Vector{a}$ is a normal vector to the hyperplane $H$ (and hence to all hyperplanes parallel to $H$).
\medskip

A hyperplane $H$ defines a \struc{positive} and a \struc{negative halfspace}

\begin{align*}
\struc{H^+} \ := \ \left\{\xb \in \R^n \suchthat \langle \Vector{a},\xb \rangle \geq b\right\} \ \text{ and } \ \struc{H^-} \ := \ \left\{\xb \in \R^n \suchthat \langle \Vector{a},\xb \rangle \leq b\right\}.
\end{align*}
\smallskip

\noindent For a given convex set $A$, $H$ is called a \struc{supporting hyperplane} if $A\cap H \neq \emptyset$ and $A \subseteq H^+$.
\medskip

Every intersection of (possibly infinitely many) positive halfspaces yields a convex set.
And, vice versa, every convex set is given by the intersection of the positive half-spaces of all its supporting hyperplanes.
Moreover, if a convex set is representable as the intersection of finitely many halfspaces, it is called a \struc{polyhedron}, and if it is moreover bounded it is called a \struc{polytope}.
\medskip

\begin{figure}[t]
  \centering
  \includegraphics[width=.4\textwidth]{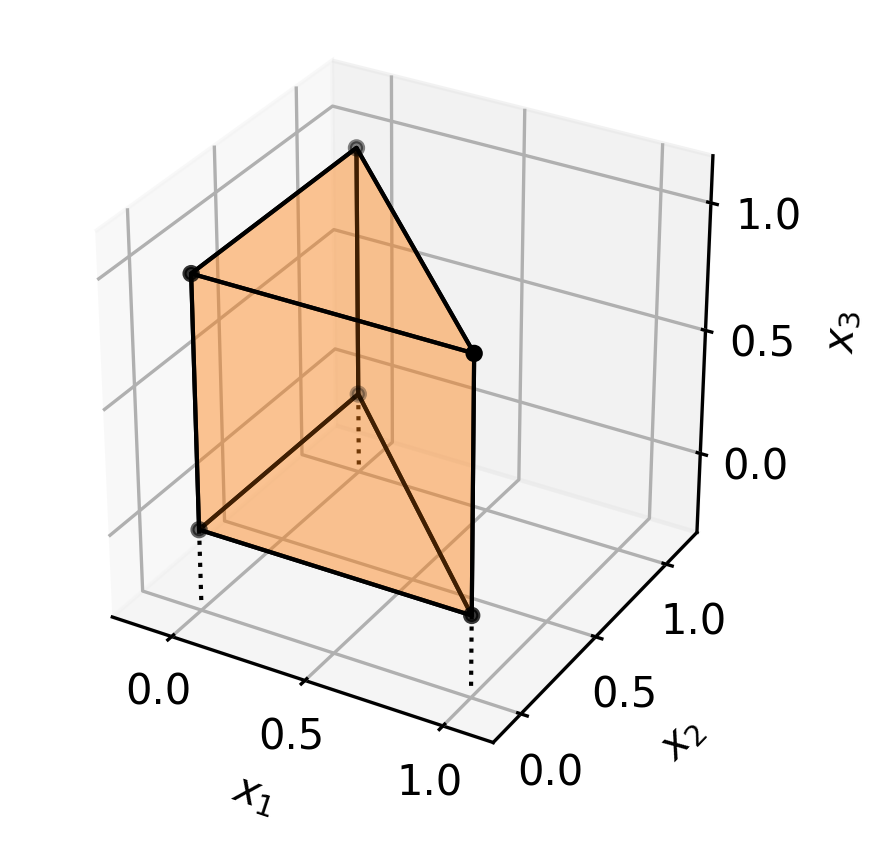}
  \caption{Three-dimensional prism over a triangle}
  \label{fig:prism}
\end{figure}

\begin{example}
	Consider the three-dimensional prism over a triangle shown in Figure~\ref{fig:prism}.
    The volume illustrated is defined by the inequalities $x_1\geq 0$, $x_2\geq 0$, $x_3\geq 0$, $x_3\leq 1$ and $x_1 + x_2 \leq 1$.
    Each inequality can be represented in terms of a positive halfspace.
    For example,
    \begin{equation*}
        \left\{\xb \in\R^3 \ \middle| \ x_1 \geq 0\right\} = \left\{\xb\in\R^3 \ \middle| \ \langle \left(\begin{smallmatrix}1\\0\\0\end{smallmatrix}\right), \xb\rangle \geq 0\right\} \ .
    \end{equation*}
    \smallskip
    
    The prism is, by definition, the set of points that satisfy all five inequalities or, in other words, the intersection of the five corresponding positive halfspaces.
	\label{Example:FirstPolytope}
\end{example}
\medskip

Polyhedra carry a significant amount of combinatorial properties and data.
For example, let $P \subset \R^n$ be a polyhedron.
As subset $F$ of $P$ is called a \struc{face} of $P$ if there exists a supporting hyperplane $H$ such that $F = P \cap H$. 
One can show that every face of a polyhedron is a polyhedron itself.
The dimension of $F$, short $\struc{\dim(F)}$ is given as the dimension of the smallest affine subspace of $\R^n$ that contains $F$.
If $P$ itself is of dimension $n$, then faces of dimension $n-1$ are called \struc{facets}, faces of dimension $1$ are called \struc{edges}, and faces of dimension $0$ are called \struc{vertices}.
\medskip

In several contexts it is useful to consider assembles of polyhedra, which intersect in a structure preserving way.
Let us make this precise.
Let $P_1,\ldots,P_k \subset \R^n$ be a collection of polyhedra.
Then the collection $\cC := \{P_1,\ldots,P_k\}$ is called a \struc{polyhedral complex} if it satisfies the following conditions:
\begin{enumerate}
	\item $\emptyset \in \cC$.
	\item If $P_i \in \cC$, then all faces of $P_i$ are also contained in $\cC$.
	\item If $P_i,P_j \in \cC$, then $P_i \cap P_j$ is a face of both $P_i$ and of $P_j$.
\end{enumerate}

\section{Classic Voronoi Diagrams}
\label{appendix:voronoi_diagrams}

This section of the appendix provides additional information about classic Voronoi diagrams.
Specifically, we informally discuss an algorithm for constructing classic Voronoi diagrams of arbitrary point sets.
First, recall Definition~\ref{def:voronoi_diagram_2d} from above.
\medskip

\begin{proposition}
    \label{proposition:voronoi_diagrams_are_polyhedral_complexes}
	Voronoi diagrams are polyhedral complexes.
	Moreover, they can be computed effectively in runtime $k \log(k)$ for point sets of cardinality $k$.
\end{proposition}
\medskip

There are different algorithms for computing Voronoi diagrams.
Below, we briefly review an algorithm that uses the \struc{Delaunay triangulation}, which is dual to the Voronoi diagram of the same set of points.
More precisely, the Voronoi diagram and Delaunay triangulation graphs of the same points are dual to each other.
As a consequence, given the Delaunay triangulation of $\pumpset$, the corresponding Voronoi diagram can be derived straightforwardly.
The graph representation of a Voronoi diagram consists of the edges and vertices of the Voronoi cells.
The Delaunay triangulation graph can be obtained by connecting any two points in $\pumpset$ with a straight line if they have the following property:
There exists a circle that passes through both points and contains no other points from $\pumpset$.
Therefore, the points in $\pumpset$ are the vertices of the Delaunay triangulation graph, and the edges are the connecting lines just defined.
\medskip

\begin{figure}[t]
  \centering
  \includegraphics[width=.45\textwidth, trim=12 12 12 12, clip]{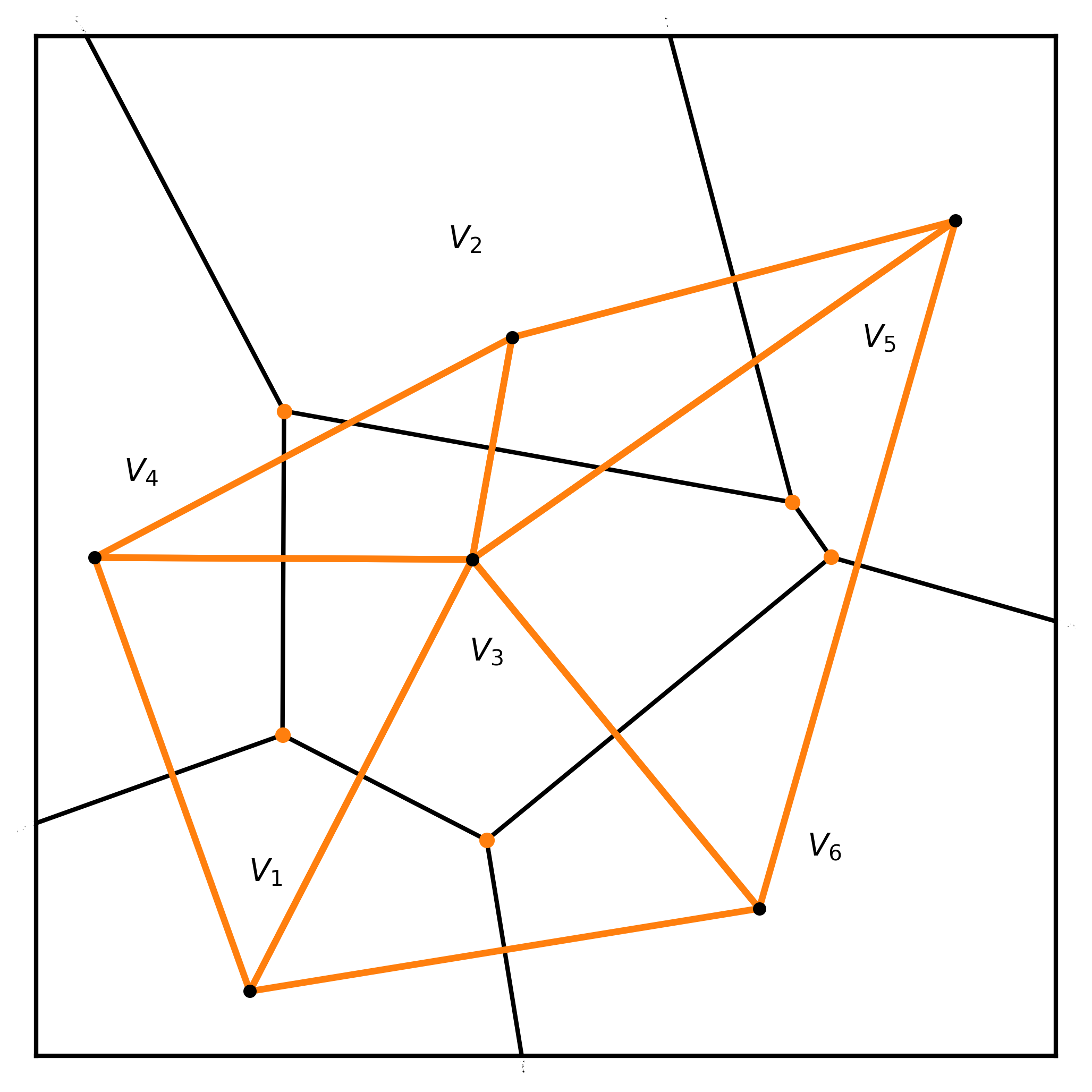}
  \caption{Voronoi diagram and Delaunay triangulation of a random set of points}
  \label{fig:classic_voronoi_diagram_random_example}
\end{figure}

\begin{example}
    We consider the Voronoi diagram from Example~\ref{example:voronoi_diagram} again.
    Figure~\ref{fig:classic_voronoi_diagram_random_example} shows the same classic Voronoi diagram (black) together with the corresponding Delaunay triangulation (orange).
\end{example}
\medskip

Before explaining how to recover a Voronoi diagram from the associated Delaunay triangulation, we describe a common algorithm for computing the Delaunay triangulation of a two-dimensional point collection.
Two dimensions are sufficient for our use case.
However, the same algorithm can be extended to higher dimensions in principle.
To ensure a unique Delaunay triangulation, we must assume that $\pumpset$ is in general position, i.e., that no three points in $\pumpset$ are co-linear and no four points in $\pumpset$ are co-circular.
\medskip

To compute the Delaunay triangulation, we first transform all points in $\pumpset$ into three dimensions, i.e., for each point $\pb=(p_1, p_2)$, we compute $\pb'=(p_1, p_2, p_1^2 + p_2^2)$.
Next, we compute the convex hull $\conv(\pumpset')$ of the associated three-dimensional point collection $\pumpset'$.
Given our assumption that $\pumpset$ is in general position, all faces of $\conv(\pumpset')$ are triangles.
Now, for each face of the convex hull, we consider an outer surface normal vector.
That is, a vector $\nb\in\R^3$ that is perpendicular to the face, points outward from the convex hull and is unique up to scaling.
If $\nb$ points downward, i.e., if $n_3$ is negative, then the face corresponds to a triangle in the Delaunay triangulation which we obtain by projecting the three vertices of the triangle back to two dimensions.
In other words, we map the vertices, which are points in $\pumpset'$, back to the corresponding points in $\pumpset$.
\medskip

The final step is to deduce the associated Voronoi diagram from the Delaunay triangulation.
First, the circumcenters of the Delaunay triangles are exactly the vertices of the Voronoi diagram.
Figure~\ref{fig:classic_voronoi_diagram_random_example} illustrates the circumcenters as orange points.
Note that depending on the type of triangle, the circumcenter may also lie outside the triangle.
To construct the edges of the Voronoi diagram, we distinguish between interior and boundary edges of the Delaunay triangulation.
Interior edges are shared by two adjacent triangles, while boundary edges belong to only one triangle.
Connecting the circumcenters of all triangles that share an interior edge yields all finite Voronoi edges (see the lines connecting the orange points in Figure~\ref{fig:classic_voronoi_diagram_random_example}). Additionally, each boundary edge of the Delaunay triangulation corresponds to an infinite Voronoi edge (a ray, to be exact).
Each infinite edge starts at the circumcenter of the respective Delaunay triangle and points outward of the Delaunay triangulation perpendicular to the respective boundary edge.
This completes the computation of the Voronoi diagram graph.
\medskip

As Proposition~\ref{proposition:voronoi_diagrams_are_polyhedral_complexes} states, Voronoi diagrams are polyhedral complexes.
In other words, each Voronoi cell $\voronoicell_i$ is a polyhedron that can be represented as the intersection of a finite number of positive halfspaces.
We need this halfspace representation for our numerical implementation, so we briefly discuss how to deduce it from the Delaunay triangulation as well.
We consider an arbitrary point $\pb_i \in \pumpset$, which is by definition also a vertex of the Delaunay triangulation.
Hence, $\pb_i$ is contained within a number of Delaunay edges and, as just discussed,
each Delaunay edge containing $\pb_i$ corresponds to an edge of the Voronoi cell $\voronoicell_i$.
Consider any point $\pb_j$ connected to $\pb_i$ by a Delaunay edge.
The \struc{bisector} of $\pb_i$ and $\pb_j$ is the set of points equidistant from $\pb_i$ and $\pb_j$.
Formally, we represent the bisector as $\{\yb\in\R^2 \mid \Vert\yb - \pb_i\Vert = \Vert\yb - \pb_j\Vert\}$.
At the same time, the bisector is a hyperplane perpendicular to the line connecting $\pb_i$ and $\pb_j$.
In other words, the vector $\pb_i - \pb_j$ is a normal vector of the bisector and the bisector has a hyperplane representation $\{\yb\in\R^2 \mid \langle\pb_i - \pb_j,\yb\rangle=b\}$.
The Voronoi cell $V_i$ can be represented as the intersection of positive halfspaces corresponding to bisectors of $\pb_i$ and adjacent points in the Delaunay triangulation.
We define

\begin{equation*}
    \ab_{ij} \coloneqq \frac{\pb_i - \pb_j}{\Vert\pb_i - \pb_j\Vert}\quad\text{and}\quad b_{ij} \coloneqq \frac{\langle \ab_{ij}, \pb_i + \pb_j\rangle}{2}
\end{equation*}
\smallskip

\noindent and therewith

\begin{equation*}
    H_{ij}^+ \coloneqq \{\yb\in\R^2 \mid \langle\ab_{ij}, \yb\rangle \geq b_{ij}\} \ .
\end{equation*}
\smallskip

\noindent The positive halfspace $H_{ij}^+$ originates from the bisector of $\pb_i$ and $\pb_j$ and covers the side of the bisector on which $\pb_i$ lies.
The only deviation from the above motivation is that the normal vector $\ab_{ij}$ now has norm one, and of course the order of $i$ and $j$ matters.
Finally, let $\neighbor_i$ denote the index set such that $j\in\neighbor_i$ if and only if there exists a Delaunay edge connecting $\pb_i$ and $\pb_j$, then

\begin{equation*}
    \voronoicell_i = \bigcap_{j\in\neighbor_i} H^+_{ij}
\end{equation*}
\smallskip

\noindent is the aimed representation of $\voronoicell_i$ as an intersection of positive halfspaces.

\section{Refined Voronoi diagrams}
\label{appendix:refined_voronoi_diagrams}

In the context of our leakage localization problem, the points determining the Voronoi cells are the locations of the vacuum connections.
The idea behind Voronoi predictors is that a leak should be closest to the particular sensor that captured the largest flow.
This concept can be refined in case there are more than two sensors.
Namely, it is likewise natural to assume that a leak should be in second closest proximity to the sensor that captured the second largest flow, and so on.

\begin{definition}
    Let $k > 2$ and $d \leq k$. Further, let
    \begin{equation*}
        \sequence_{d} \coloneqq \{\tb = (i_{1},\dots, i_{d}) \mid i_{\alpha} \in [k]  \text{ and } i_{\alpha}\neq i_{\beta} \text{ for } \alpha\neq \beta\}
    \end{equation*}
    be the set of length $d$ tuples from $[k]$ without repetition. For any $\tb\in \sequence_{d}$, let
    \begin{equation*}
        \setsequence(\tb) \coloneqq \{i_{1},\dots,i_{d}\}\quad\text{and}\quad\abs{\tb}\coloneqq\abs{\setsequence(\tb)}
    \end{equation*}
    be the set of elements in $\tb$ and its cardinality. Then, we define the corresponding \struc{refined Voronoi cell} as
    \begin{equation*}
        \label{eq:refined_voronoi_region}
        \struc{\voronoicell_{\tb}} \coloneqq \{\yb \in \R^2 \mid \Vert\yb - \pb_{i_{1}}\Vert\leq \dots \leq \Vert\yb - \pb_{i_{\abs{\tb}}}\Vert\leq \Vert\yb - \pb_{j}\Vert \text{ for all } j \in [k] \setminus \setsequence(\tb)\} .
    \end{equation*}
    We call the collection of all refined Voronoi cells based on $d$-tuples the \struc{refined Voronoi diagram of order $d$} of $\pumpset$ and denote it by $\struc{\generalizedvoronoidiagram_d(\pumpset)} = (\voronoicell_{\tb})_{\tb \in \sequence_{d}}$.
\end{definition}

\begin{figure}[t]
  \centering
  \includegraphics[width=.45\textwidth, trim=12 12 12 12, clip]{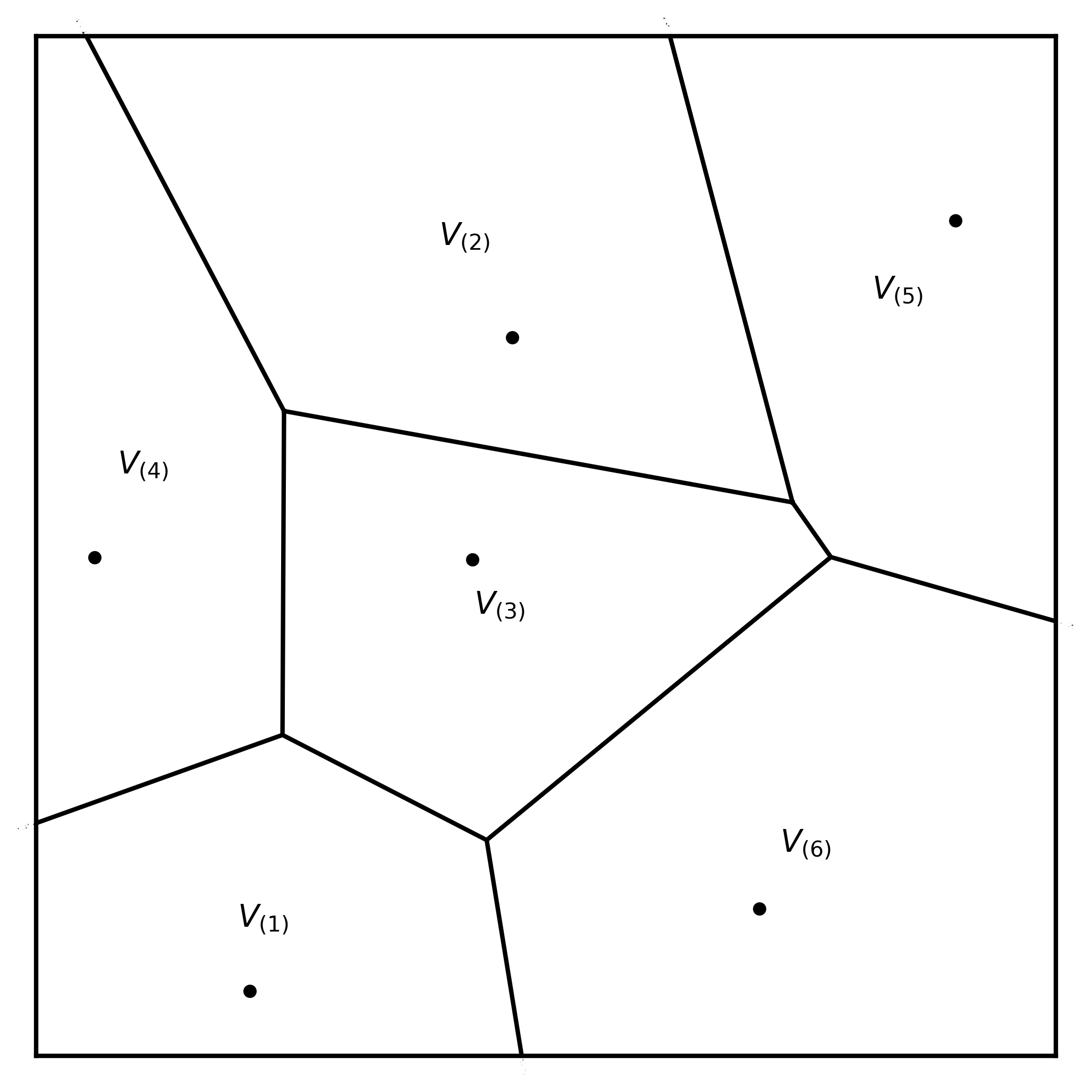}\quad
  \includegraphics[width=.45\textwidth, trim=12 12 12 12, clip]{voronoi_permutation_example.png}
  \caption{Classic Voronoi diagram of a random set of points / Refined Voronoi diagram of the same points}
  \label{fig:refined_voronoi_diagram_random_example}
\end{figure}

\begin{example}
    $\voronoicell_{\tb}$ is the region of points that are closest to $\pb_{i_1}$, second closest to $\pb_{i_2}$, and so on.
    Note that $\generalizedvoronoidiagram_1(\pumpset) = (\voronoicell_{(1)}, \dots, \voronoicell_{(k)})$ is exactly the classic Voronoi diagram $\voronoidiagram(\pumpset) = (\voronoicell_1, \dots, \voronoicell_k)$, as can be seen by comparison of Figures~\ref{fig:classic_voronoi_diagram_random_example} and \ref{fig:refined_voronoi_diagram_random_example} (left).
    The corresponding refined Voronoi diagram of order two $\generalizedvoronoidiagram_2(\pumpset)$ is illustrated on the right side of Figure~\ref{fig:refined_voronoi_diagram_random_example}.
    The case $d=k$ corresponds to a maximal refinement of Voronoi regions by means of $k$-tuples $\tb\in \sequence_{k}$.
\end{example}
\medskip

In the following, we assume that an algorithm to compute classic Voronoi diagrams is available.
Based on that, we derive a scheme that can be used to compute refined Voronoi diagrams recursively.
The proposed scheme is outlined in Algorithm~\ref{alg1}.
Assuming that $\generalizedvoronoidiagram_d(\pumpset)$ for $d\leq k-1$ is given as input, the refined Voronoi diagram of subsequent order $\generalizedvoronoidiagram_{d+1}(\pumpset)$ results as output.
\medskip

Each tuple $(i_{1},\dots,i_{d})$ is visited during the outer for loop, and the inner for loop then provides all refined Voronoi cells for tuples $(i_{1},\dots,i_{d},j)$ with $j\in [K] \setminus \{i_{1},\dots,i_{d}\}$.
To that end, a standard Voronoi diagram with reduced point set

\begin{equation*}
    \voronoidiagram(\pumpset\setminus\{\pb_{i_1},\dots,\pb_{i_d}\}) \eqqcolon (\voronoicellreduced_i)_{i\in[k]\setminus\{i_1,\dots,i_d\}}
\end{equation*}
\smallskip

\noindent is computed first.
Each new Voronoi region of subsequent order is then obtained by taking the intersection

\begin{equation*}
    \voronoicell_{(i_{1},\dots,i_{d},j)} \coloneqq \voronoicell_{(i_{1},\dots,i_{d})} \cap \voronoicellreduced_{j} \ .
\end{equation*}
\smallskip

\noindent The interpretation of that step is as follows:
We have previously computed the cell $\voronoicell_{(i_{1},\dots,i_{d})}$.
All points in that cell are closest to $\pb_{i_{1}}$, second closest $\pb_{i_{2}}$, and so on.
In the current step, we are looking for the $k - d$ subsets of $\voronoicell_{(i_{1},\dots,i_{d})}$ that are $d+1$-st closest to $\pb_{j}$ for $j\in [k]\setminus \{i_{1},\dots,i_{d}\}$.
These are obtained through the intersections with $\voronoicellreduced_{j}$.

\begin{algorithm}[t]
  \caption{Computation of $(\voronoicell_{\tb})_{\tb\in \sequence_{d+1}}$ given $(\voronoicell_{\tb})_{\tb\in \sequence_{d}}$}
  \label{alg1}
  \begin{algorithmic}
    \ForAll{$(i_{1},\dots,i_{d}) \in \sequence_{d}$}
    \State Compute Voronoi diagram $(\voronoicellreduced_{i})_{i\in [k] \setminus \{i_{1},\dots,i_{d}\}}$
    \ForAll{$j\in [k]\setminus \{i_{1},\dots,i_{d}\}$}
    \State $\voronoicell_{(i_{1},\dots,i_{d}, j)} \coloneqq \voronoicell_{(i_{1},\dots,i_{d})} \cap \voronoicellreduced_{j}$
    \EndFor
    \EndFor
  \end{algorithmic}
\end{algorithm}

\begin{figure}[t]
    \centering
    \includegraphics[scale=1.1, trim=8 8 8 8, clip]{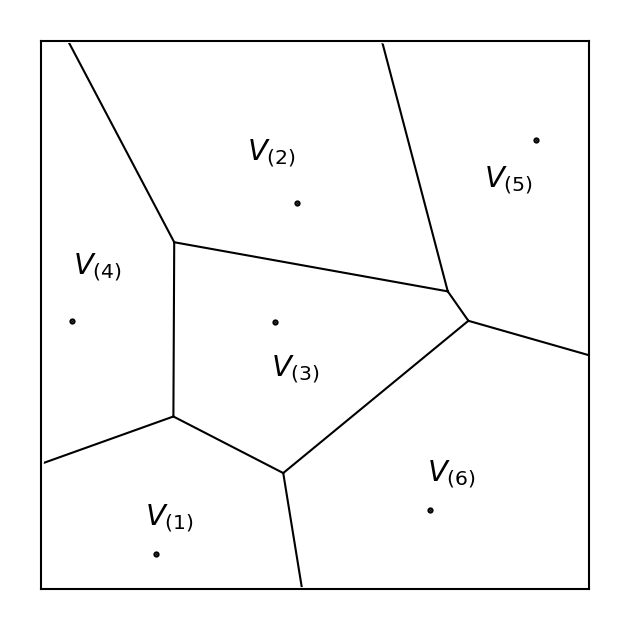}
    \includegraphics[scale=1.1, trim=8 8 8 8, clip]{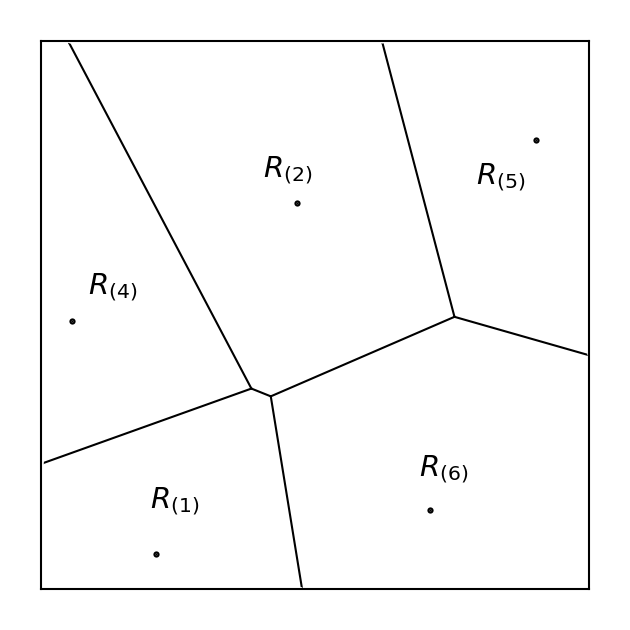}
    \includegraphics[scale=1.1, trim=8 8 8 8, clip]{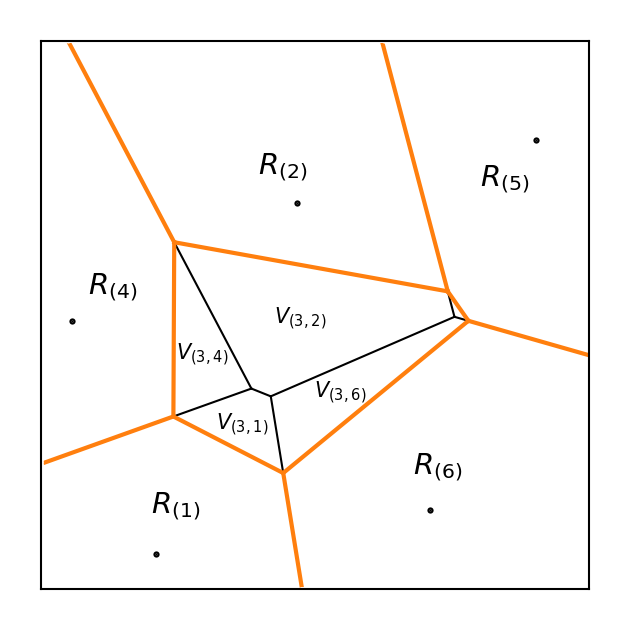}
    \caption{Subsequent steps of Algorithm~\ref{alg1} applied to a random example.}
    \label{fig:voronoi_generalization_algorithm_steps}
\end{figure}

\begin{example}
  Figure~\ref{fig:voronoi_generalization_algorithm_steps} illustrates the transition from $d=1$ to $d=2$ using the example from Figure~\ref{fig:refined_voronoi_diagram_random_example}.
  The cell $\voronoicell_{(3)}$ is subdivided into $\voronoicell_{(3, 1)}$, $\voronoicell_{(3, 2)}$, $\voronoicell_{(3, 4)}$, $\voronoicell_{(3, 5)}$ and $\voronoicell_{(3, 6)}$.
\end{example}

As discussed above, classic Voronoi diagrams are polyhedral complexes, and their cells can be represented as the intersection of a finite number of positive halfspaces.
These statements also hold for refined Voronoi diagrams.
To understand how a refined Voronoi cell can be represented as the intersection of halfspaces in the case of $d=2$, consider the definition according to the inner loop of Algorithm~\ref{alg1}:

\begin{equation*}
    \voronoicell_{(i, j)} \coloneqq \voronoicell_{i} \cap \voronoicellreduced_{j}
\end{equation*}
\smallskip

\noindent In other words $\voronoicell_{(i, j)}$ is the intersection of $\voronoicell_i$ and $\voronoicellreduced_j$, both of which are classic Voronoi cells, the latter only with respect to a reduced point set.
Let $\neighbor_j'$ include the indices of the neighbors of $\pb_j$ in the Delaunay triangulation of the reduced point set. Then, analogous to above, the two cells can be represented as

\begin{equation*}
    \voronoicell_i = \bigcap_{\mu\in\neighbor_i} H^+_{i\mu}\quad\text{and}\quad \voronoicellreduced_j = \bigcap_{\nu\in\neighbor_j'} H^+_{j\nu}
\end{equation*}
\smallskip

\noindent and hence,

\begin{equation*}
    \voronoicell_{(i,j)} = \left(\bigcap_{\mu\in\neighbor_i} H^+_{i\mu}\right) \ \bigcap \ \left(\bigcap_{\nu\in\neighbor_j'} H^+_{j\nu}\right) \ .
\end{equation*}
\smallskip

Thus, the set of halfspaces to represent $\voronoicell_{(i, j)}$ is simply the union of the respective sets for $\voronoicell_i$ and $\voronoicellreduced_j$.
\bigskip

\end{document}